\documentclass[10pt]{article}
\usepackage{amscd,amsfonts,amssymb,amscd,amsmath,latexsym}
\usepackage[hypertex,backref]{hyperref}
\usepackage[all]{xy}
\makeatletter
\renewcommand{\subsubsection}{\@startsection
{subsubsection}
{1}
{0mm}
{0mm}
{0mm}
{\normalfont\normalsize\itshape}}
\makeatother
  
\usepackage{mathrsfs}
\newtheorem{theorem}{Theorem}[section] 
\newtheorem{prop}[theorem]{Proposition}
\newtheorem{lem}[theorem]{Lemma}
\newtheorem{ddd}[theorem]{Definition}

\newcommand{\Mor}{\mathrm{Mor}}

\newcommand{\colim}{\mathrm{colim}}

\newcommand{\forget}[1]{}

{\catcode`@=11\global\let\c@equation=\c@theorem}

\newcommand{\Z}{\mathbb{Z}}

\newcommand{\proof}{{\it Proof.$\:\:\:\:$}}

\newcommand{\R}{\mathbb{R}}

\newcommand{\Q}{\mathbb{Q}}

\newcommand{\bL}{{\bf L}}

\newcommand{\bN}{\mathbf{N}}

\newcommand{\cH}{\mathcal{H}}

\newcommand{\cE}{\mathcal{E}}

\newcommand{\cC}{\mathcal{C}}
\newcommand{\cK}{\mathcal{K}}
\newcommand{\cA}{\mathcal{A}}

\newcommand{\Hom}{{\tt Hom}}

\newcommand{\id}{{\tt id}}

\def\imath{{i}}

\def\hB{\hspace*{\fill}$\Box$ \newline\noindent}

\newcommand{\cS}{\mathcal{S}}

\def\hB{\hspace*{\fill}$\Box$ \\[0cm]\noindent}

\newcommand{\bH}{\mathbf{H}}

\newcommand{\pr}{{\tt pr}}
\newcommand{\bX}{\mathbf{X}}
\newcommand{\bY}{\mathbf{Y}}

\newcommand{\ev}{{\tt ev}}

\newcommand{\Sets}{\mathcal{S}ets}

\title{Sheaf theory for stacks in manifolds and twisted cohomology for $S^1$-gerbes}
\author{Ulrich Bunke, Thomas Schick and Markus Spitzweck
\thanks{Mathematisches Institut, Universit{\"a}t G{\"o}ttingen,
Bunsenstr. 3-5, 37073 G{\"o}ttingen, GERMANY,  bunke@uni-math.gwdg.de,
schick@uni-math.gwdg.de, spitz@uni-mah.gwdg.de}
}

\begin{document}
\newcommand{\Top}{{\tt Top}}
\newcommand{\bM}{{\mathbf{M}}}
\newcommand{\ori}{\tt or} 
\newcommand{\cat}{{\tt cat}}
\newcommand{\cov}{{\tt cov}}
\newcommand{\Ab}{{\tt Ab}}
\newcommand{\Sh}{{\tt Sh}}
\newcommand{\Site}{{\tt Site}}
\newcommand{\Mf}{{\tt Mf}}
%\begin{Abstract}
%\end{Abstract}

\maketitle 
\tableofcontents

\newcommand{\bG}{\mathbf{G}}

\section{Introduction}

\subsection{About the motivation}

\subsubsection{}

Given a closed three form $\lambda\in \Omega^3(X)$ on a smooth manifold $X$, the usual definition of twisted de Rham cohomology is as the cohomology of the two-periodic complex $(\Omega^\bullet_{per}(X),d_\lambda)$, where $$\Omega^\bullet_{per}(X):=\bigoplus_{n\in \Z}\Omega^{\bullet+2n}(X)\ ,\quad \mbox{and}\quad  d_\lambda:=d_{dR}+\lambda$$ is the sum of the de Rham differential and the multiplication operator by the form $\lambda$.

\subsubsection{}

Twisted de Rham cohomology is in particular interesting as a target of the Chern character from twisted $K$-theory. In this case 
$[\lambda]\in H^3(X;\R)$ is the real image of an integral class $\lambda_\Z(P)\in H^3(X;\Z)$ which classifies a principal bundle $P\to X$ with structure group $PU$, the projective unitary group of a complex infinite-dimensional separable Hilbert space. The twisted $K$-theory depends functorially on $P$ in a non-trivial  manner.

The twisted cohomology as defined above depends on the cohomology class $[\lambda]$ up to (in general)  non-canonical isomorphism.
The draw-back of  this definition of twisted cohomology above is that it is not functorial in the twist $P\to X$ of $K$-theory since there no canonical choice of a three-form $\lambda$ representing the image of $\lambda_\Z(P)$ in real cohomology.

\subsubsection{}
 The main goal of the present note is to propose an alternative functorial definition of the twisted cohomology as the real cohomology of a stack $G_P$ which is canonically associated to the $PU$-bundle $P\to X$. The stack $G_P$ is the stack of $U$-liftings of $P\to X$, where $U$ is the unitary group of the Hilbert space and $U\to PU$ is the canonical projection map. It is also called the lifting gerbe of $P$.

In order to define the cohomology of a stack like $G_P$ we develop a sheaf theory set-up for stacks in smooth manifolds. Our main result 
Theorem \ref{main} is the key step in the verification that the cohomology according to the new sheaf-theoretic definition is essentially isomorphic (non-canonically) to the twisted cohomology as defined above.

We have chosen to work with stacks in smooth manifolds since  we are heading towards a comparison with de Rham cohomology. A parallel theory
can be set up in the topological context. Together with applications to
$T$-duality and delocalized cohomology it will  be discussed in detail in the
subsequent papers %\cite{bssf},
 \cite{bss1} and \cite{bssm}.

\subsubsection{}

In \cite{MR2172499, math.DG/0605694}, a different version of sheaf theory and
cohomology of 
stacks is developed. Already the site associated to a stack in these papers is different from ours, as we will discuss later (compare \ref{site-comp}). But, there is a
comparison map which in the situations we are interested in (in particular for
constant sheaves and the de Rham sheaf) induces an isomorphism in cohomology.

We have to develop our own version of sheaf theory and sheaf cohomology for
stacks, because our argument heavily relies on functorial constructions associated to maps between stacks.
This calculus has not been developed in the references above.
% properties ---which partly we have to develop in an ad hoc fashion which works
% precisely in the given context. 
%These kind of constructions is not carried out and analyzed in the work of Behrend et al.

\subsubsection{}

The twists for our new cohomology theory are smooth gerbes $G\to X$ with band $U(1)$. The lifting gerbe $G_P\to X$ of a $PU$-bundle mentioned above is an example. Advantages of our new definition are:
\begin{enumerate}
\item The twisted cohomology depends  functorially on the twist.
\item One can define twisted cohomology with coefficients in an arbitrary abelian group.
\item The definition can easily be generalized to the topological context.
\end{enumerate}

\subsubsection{}

In Subsection \ref{rtrt1} we give a complete technical statement of our main result written for a reader familar with  the language of stacks, sites, and sheaf theory.
The third  part of the introduction, Subsection \ref{se3}, is devoted to a detailed motivation with references to the literature and a less technical introduction of the language and the description of the result. Finally, Subsection \ref{introsch} is an introduction to the technical sheaf theoretic part of the present paper.

\subsection{Statement of the main result}\label{rtrt1}

\subsubsection{}

We consider a stack $G$ on the category of smooth manifolds equipped with the usual topology of open coverings.
To $G$ we associate a site $\bG$ as a subcategory of manifolds over $G$. The objects of this site are representable smooth maps
$U\rightarrow G$ from smooth manifolds to $G$. 
A covering $(U_i\rightarrow U)_{i\in I}$  is a collection of morphisms which are submersions and such that
$\sqcup_{i\in I}U_i\rightarrow U$ is surjective (see \ref{site33} for a precise definition).

\subsubsection{}

To the site $\bG$ we associate
the categories of presheaves $\Pr\bG$ and sheaves $\Sh\bG$ of sets as well as the lower bounded derived categories  $D^+(\Pr_{\Ab} \bG)$ and $D^+(\Sh_{\Ab}\bG)$ of the abelian categories $\Pr_{\Ab}\bG$ and $\Sh_{\Ab}\bG$ of presheaves and sheaves of abelian groups.

\subsubsection{}

Let $i:\Sh\bG\rightarrow \Pr\bG$ be the natural inclusion, and let $i^\sharp:\Pr\bG\rightarrow \Sh\bG$ be its left adjoint, the sheafification functor. As a right adjoint the functor $i$ is left exact and admits a right derived functor $Ri:D^+(\Sh_\Ab\bG)\to D^+(\Pr_\Ab\bG)$.

\subsubsection{}
If $G\rightarrow X$ is a morphism of stacks, then we define a functor
$f_*:\Pr\bG\rightarrow \Pr\bX$. Note that if $f$ is not representable, then this map is not associated to a map of sites. If $F\in \Pr\bG$ and $(U\rightarrow X)\in \bX$, then
we set  (see \ref{dd3})
$$f_*(U):=\lim F(V)\ ,$$
where the limit is taken over the category of diagrams
$$\xymatrix{V\ar[r]\ar[d]&G\ar[d]^f\\U\ar@{=>}[ur]\ar[r]&X}\ .$$
It turns out that $f_*$ admits a left adjoint. Therefore it is left exact and  admits a right derived functor
$Rf_*:D^+(\Pr_{\Ab}\bG)\rightarrow D^+(\Pr_{\Ab}\bX)$.

\subsubsection{}\label{oldes}

Let
$f:G\rightarrow X$ be a smooth gerbe with band $S^1$ over the smooth manifold $X$.
We consider the sheafification $i^\sharp\R_\bG$ of the constant
presheaf $\R_\bG$ on $\bG$ with value $\R$. Our main result describes
$$i^\sharp \circ Rf_*\circ Ri(i^\sharp \R_\bG)\in D^+(\Sh_{\Ab}\bX)$$ in terms of a deformation of the de Rham complex.

The gerbe $f:G\rightarrow X$ is classified by a Dixmier-Douday class $\lambda_\Z\in H^3(X;\Z)$.
Let $\lambda\in \Omega^3(X)$ be a closed form such that $ [\lambda]\in H^3(X;\R)$ represents the image of $\lambda_\Z$ under $H^3(X;\Z)\rightarrow H^3(X;\R)$.

For a manifold $X$ the objects $(U,p)$ of the site $\bX$ are submersions $p:U\rightarrow X$ from smooth manifolds $U$ to $X$. This differs from the usual convention, where the site is the category of open subsets of $X$.

We form the complex of presheaves $(U,p)\mapsto \Omega^\cdot[[z]]_\lambda(U,p)$ on $\bX$, which associates to $(U,p)\in \bX$ the complex of formal power series of smooth real differential forms on $U$ with differential 
$$d_\lambda:=d_{dR}+T\lambda\ ,$$ where $z$ is a formal variable of degree $2$, $T:=\frac{d}{dz}$,
$d_{dR}$ is the de Rham differential, and $\lambda$ stands for multiplication by $p^*\lambda$. It turns out that this is actually a complex of sheaves (see Lemma \ref{resd97}).

\subsubsection{}

The main result of the present paper is the following theorem.
\begin{theorem}\label{main}
In $D^+(\Sh_{\Ab}\bX)$ we have an isomorphism $i^\sharp \circ Rf_*\circ Ri(i^\sharp \R_\bG)\cong \Omega^\cdot[[z]]_\lambda$.
\end{theorem}

\subsubsection{}

The projection map $f:G\rightarrow X$ of a gerbe is not representable so that
$f_*:\Pr\bG\rightarrow \Pr\bX$ does not come from an associated map of sites.
Therefore, in order to define $Rf_*$ and to verify the theorem we have to develop some standard elements of sheaf theory for stacks in smooth manifolds. This is the contents of Section \ref{se1} (see \ref{introsch} for an introduction). 
In Section \ref{se2} we verify Theorem \ref{main}.

\subsection{Twisted cohomology and gerbes}\label{se3}

\subsubsection{}

A closed  three-form  $\lambda\in \Omega^3(X)$ on a smooth manifold $X$ can be used to perturb the de Rham differential $$d_{dR}\:\:\:\\\leadsto\:\:\: d_{dR}+\lambda=:d_\lambda\ .$$ The cohomology of the two-periodic complex 
$$\dots\stackrel{d_\lambda}{\rightarrow} \Omega^{even}(X) \stackrel{d_\lambda}{\rightarrow} \Omega^{odd}(X)\stackrel{d_\lambda}{\rightarrow} \Omega^{even}(X)\stackrel{d_\lambda}{\rightarrow}\dots$$
is called the $\lambda$-twisted cohomology of $X$ and often denoted by $H^*(X;\lambda)$.
This ad-hoc definition appears in various places in the recent mathematical literature
(let us mention just  \cite{math.KT/0510674},  \cite{MR1911247}, \cite{MR1977885}, \cite{MR2080959}) and in the physics literature. 
A closely related and essentially equivalent definition
\cite{math.KT/0404329} uses the complex $(\Omega^\cdot(X)((u)),d_{dR}-u\lambda)$, where $u$ is a formal variable of degree $-2$, and $''((u))''$ stands for formal Laurent series.

\subsubsection{}

It is known that
the isomorphism class of the $\lambda$-twisted cohomology group only depends on the cohomology class $[\lambda]\in H^3(X;\R)$. 
If $f:Y\rightarrow X$ is a smooth map, then we have a functorial map
$f^*:H^*(X;\lambda)\rightarrow H^*(Y;f^*\lambda)$ which  essentially only depends on the homotopy class of $f$. Furthermore, $\lambda$-twisted cohomology has a Mayer-Vietoris sequence and is a module over $H^*(X;\R)$. It now appears as a natural question to understand $\lambda$-twisted cohomology as a concept  of  algebraic topology. 

\subsubsection{}

One attempt is the approach
of \cite{math.AT/0206257} in which the complex of smooth differential forms is replaced by similar objects in algebraic topology. 

The proposal of \cite{math.KT/0510674} to use the singular de Rham complex goes into the same direction. Observe that we can use the filtration of $\Omega^{ev}(X)$ and $\Omega^{odd}(X)$ by degree in order to construct a spectral sequence converging to $H^*(X;\lambda)$. Its $E_2$-page  involves $H^{*}(X;\R)$ (as $\Z/2\Z$-graded vector spaces). The next possibly non-trivial differential of this spectral sequence is the multiplication by the class $[\lambda]$. In \cite{math.KT/0510674} the higher differentials of this spectral sequence are identified as Massey products.

\subsubsection{}

A natural homotopy theoretic framework for twisted cohomology theories
would be some version of parametrized stable homotopy theory as developed e.g. in \cite{math.AT/0411656}. In such a theory a twist of a generalized cohomology theory (represented by a spectrum $E$) is a parametrized
spectrum $\cE$ over $X$ with typical fibre equivalent to $E$ (think of a bundle of spectra). The twisted cohomology groups $H^*(X;\cE)$ are then given by the homotopy groups of the spectrum of sections of $\cE$.
In order to interpret $\lambda$-twisted cohomology in this manner one would have to relate three-forms on $X$ with parametrized versions of the Eilenberg-MacLane spectrum $H\R$.

Let us mention that alternatively to \cite{math.AT/0411656} other reasonable versions of a stable homotopy theory over $X$  could be based on presheaves of spectra over $X$ or $\Omega(X)$-equivariant spectra, where $\Omega(X)$ denotes the based loops of $X$.

\subsubsection{}

One motivation for introducing $\lambda$-twisted cohomology is based on the fact that it can be used as a target of the Chern character from twisted $K$-theory. It is known that $H^3(X;\Z)$ classifies a certain subset of isomorphisms classes of parametrized spectra $\cK$ with fibre equivalent to the complex $K$-theory spectrum $K$. This follows from the splitting $BGL_1(K)\cong K(\Z,3)\wedge T$.
Here $GL_1(K)$ denotes the grouplike monoid of units of the $K$-theory spectrum, 
$K(\Z,3)$ denotes an Eilenberg-MacLane space, and $T$ is an auxiliary space. We refer
to \cite{MR0494077} for more details. Chern characters are constructed in \cite{MR1911247}, 
\cite{math.KT/0510674}, \cite{MR1977885}, \cite{math.KT/0404329}. Note that in these works twisted $K$-theory is not defined in homotopy theoretic terms but using sections in bundles of Fredholm operators, bundle gerbe modules or $K$-theory of $C^*$-algebras.
 If $\lambda_\Z\in H^3(X;\Z)$ classifies the parametrized $K$-theory spectrum $\cK$, then
the Chern character has values in $H^*(X;\lambda)$, where $[\lambda]$ is the image of $\lambda_\Z$ under the map $H^3(X;\Z)\rightarrow H^3(X;\R)$. Such a definition can not be natural since in general $\cK$ has non-trivial automorphisms which are not reflected by $H^3(X;\lambda)$. 

A completely natural definition of a Chern character with values even in a twisted rational cohomology could be induced from the canonical rationalization map $\cK\rightarrow \cK_\Q$ if we like to define twisted rational cohomology using $\cK_\Q$.

\subsubsection{}

Above we have seen that $H^3(X;\Z)$ classifies a subset of the isomorphism classes of parametrized $K$-theory spectra over $X$. This can in fact be seen directly. Let $U$ be the unitary group of a separable infinite  dimensional complex Hilbert space. Equipped with the topology induced by the operator norm it is a topological group. By Kuiper's theorem it is contractibe so that the
projective unitary group $PU:=U/U(1)$ has the homotopy type of $BU(1)\cong K(\Z,2)$. Taking the classifying space once more we have $BPU\cong K(\Z,3)$.
This shows that $H^3(X;\Z)$ classifies isomorphism classes
of $PU$-principal bundles over $X$. One can now manufacture a $PU$-equivariant version
of a $K$-theory spectrum $K$ (see e.g. \cite{MR2122155}). If $P\rightarrow X$ is a $PU$-principal bundle, then one can define the bundle of spectra $\cK:=P\times_{PU} K$ over $X$.
Alternatively one could construct twisted $K$-theory starting from a bundle of projective Hilbert spaces as in \cite{MR2172633}. As a result of this discussion one should consider $PU$-principal bundles as more primary objects.

\subsubsection{}\label{rrtp}

The theory of bundle gerbes initiated in \cite{MR1405064} and continued  in \cite{MR1794295} aims at a categorification of $H^3(X;\Z)$ in a similar manner as $U(1)$-principal bundles categorify $H^2(X;\Z)$. The $PU$-principal bundles considered above are particularly nice examples of bundle gerbes. Other examples of bundle gerbes are introduced in \cite{MR1876068}. In order to simplify we forget the smooth structure of $X$ for the moment and work in the category of topological spaces.

Let us represent $X$ as a moduli space of a groupoid
$A^1\Rightarrow A^0\rightarrow X$  in topological spaces, i.e. we represent $X$ as the quotient of the space of objects $A^0$ by the equivalence relation $A^1$. In addition we shall assume that the range and source maps have local sections.
Then a bundle gerbe is the same as a central $U(1)$-extension
$$\xymatrix{U(1)\ar[d]&&\\ \tilde A^1\ar@{=>}[r]\ar[d]& A^0\ar@{=}[d]&\\A^1\ar@{=>}[r]&A^0\ar[r]&X}$$
of topological groupoids.

In order to relate the  $PU$-principal bundle $P\rightarrow X$  with a bundle gerbe we represent
$X$ as the moduli space of the action groupoid $P\times PU\Rightarrow P\rightarrow X$. The central  $U(1)$-central extension of this groupoid is given by $P\times U\Rightarrow P$.

\subsubsection{}

The picture of a gerbe in \cite{MR1876068} is obtained by choosing an open covering $(U_i)_{i\in I}$ of $X$ 
and forming the representation
$$\bigsqcup_{i,j}U_i\cap U_j\Rightarrow \bigsqcup_{i\in I} U_i\rightarrow X\ .$$
The data of a $U(1)$-central extension of this groupoid is equivalent to transition line bundle data and trivializations over triple intersection  considered in \cite{MR1876068}.

One can build a two-category of topological groupoids by inverting Morita equivalence such that equivalence classes of $U(1)$-central extensions of groupoids representing $X$ are indeed classified by $H^3(X;\Z)$ (see e.g. \cite{math.KT/0306138}).

\subsubsection{}

A more natural view on this category of groupoids is through stacks on topological spaces $\Top$. We consider $\Top$ as a Grothendieck site where covering families are given by coverings by families of open subspaces.
 
Note that groupoids form a two-category. A stack $G$ on $\Top$ can be viewed as an object which associates 
to each space $U\in \Top$ a groupoid $G(U)$,
to a morphism $U^\prime\rightarrow U$ a homomorphism of groupoids $G(U)\rightarrow G(U^\prime)$, to a chain of composable morphisms
$$\xymatrix{U^{\prime\prime}\ar[dr]\ar[rr]&&U\\
&U^\prime\ar[ur]&}$$
a two-isomorphism
$$\xymatrix{G(U)\ar[dr]\ar[rr]&&G(U^{\prime\prime})\\
&G(U^\prime)\ar[ur]\ar@{=>}[u]&}$$
satisfying a natural associativity relation, and such that
$G$ satisfies descent conditions for the covering families of $U$.
Precise definitions can be found e.g. in \cite{math.AG/0503247}, \cite{heinloth}, \cite{MR1197353}. A space $V\in \Top$ can be viewed as a stack by the Yoneda embedding such that $V(U)=\Hom_{\Top}(U,V)$ (where we consider sets as groupoids with only identity morphisms).

\subsubsection{}\label{lg}

As an illustration we explain a canonical construction which associates to a $PU$-principal bundle $P\rightarrow X$ over a space $X$ a stack $G_P$ together with a map $G_P\to X$. It will be called the lifting gerbe of $P$.

Observe that $U$ acts on $P$ via the canonical homomorphims $U\to PU$. 
For a space $T\in \Top$ the objects of the groupoid $G_P(T)$ are the diagrams
$$\xymatrix{Q\ar[r]\ar[d]&P\ar[d]\\
T\ar[r]&X}\ ,$$ 
where $Q\to T$ is a $U$-principal bundle, and
$Q\to P$ is $U$-equivariant.

A morphism between two such objects
$$\xymatrix{Q\ar[r]\ar[d]&P\ar[d]\\
T\ar[r]&X}\ ,\quad \xymatrix{Q^\prime \ar[r]\ar[d]&P\ar[d]\\
T\ar[r]&X}\ ,$$
 is an isomorphism of $U$-principal bundles
$Q\to Q^\prime$ over $T$ which is compatible with the maps to $P$.

Finally, for a map $T^\prime\to T$ the functor
$G_P(f):G_P(T)\to G_P(T^\prime)$ maps the object
$$\xymatrix{Q\ar[r]\ar[d]&P\ar[d]\\
T\ar[r]&X}\in G_P(T)$$ to the induced diagram  $$\xymatrix{T^\prime\times_T Q \ar[r]\ar[d]&P\ar[d]\\
T^\prime\ar[r]&X}\in G_P(T^\prime)\ ,$$ and a morphism $Q\to Q^\prime$ to the induced morphism $T^\prime\times_TQ\to T^\prime\times_TQ^\prime$.
We leave it as an exercise to check that this presheaf of groupoids is a stack.

The morphism $G_P\to X$  maps the object
$$\xymatrix{Q\ar[r]\ar[d]&P\ar[d]\\
T\ar[r]&X}\in G_P(T)$$ to the underlying map $T\to X$ which is considered as an element of $X(T)$.

\subsubsection{}

A diagram of $PU$-principal bundles 
$$\xymatrix{P\ar[r]\ar[d]&P^\prime\ar[d]\\X\ar[r]&X^\prime}$$
functorially induces a diagram of stacks
$$\xymatrix{G_P\ar[r]\ar[d]&G_{P^\prime}\ar[d]\\X\ar[r]&X^\prime}$$
in the obvious way.

\subsubsection{}

A topological groupoid $A:A^1\Rightarrow A^0$ represents a stack
$[A^1/A^0]$ on topological spaces. It associates to each space $U$ the groupoid $[A^0/A^1](U)=\Hom(U,[A^0/A^1])$ of $A$-principal bundles on $X$ and isomorphisms (see \cite{heinloth}). A morphism of groupoids gives rise via an associated bundle construction to a map of stacks. As discussed in \cite{MR1401424}
one can embed in this way
the two-category of topological groupoids (with Morita equivalence inverted) mentioned at the end of \ref{rrtp} as a full subcategory of stacks on $\Top$.
The image of this embedding consists of topological stacks $G$, i.e. stacks which admit an atlas $A^0\rightarrow G$.

An atlas is a surjective
representable morphism $A^0\rightarrow G$ admitting local sections, where $A^0$ is a space. Given an atlas of $G$ we can construct a groupoid $A^1\Rightarrow A^0$. The morphism space of the groupoid is given by  $A^1:=A^0\times_GA^0$.
We then have an equivalence of stacks  $[A^0/A^1]\cong G$.

A map of stacks $G\rightarrow H$ is called representable if for any map $U\rightarrow H$ with $U$ a space $U\times_HG$ is equivalent to a space. The representability condition on $A^0\rightarrow G$ ensures that $A^1:=A^0\times_GA^0$ is a space.

\subsubsection{}\label{atla}

The lifting gerbe $G_P$ of a $PU$-principal bundle \ref{lg} is a topological stack. In order to construct an atlas we choose a covering of $X$ by open subsets on which $P$ is trivial. Let $A$ be the disjoint union of the elements of the convering, and $A\to X$ be the canonical map. By choosing  local trivializations we obtain the lift in the diagram
$$\xymatrix{&P\ar[d]\\A\ar[r]\ar@{.>}^s[ur]^s&X}\ .$$
We now consider the diagram
$$\xymatrix{A\times U\ar[r]^\phi \ar[d]&P\ar[d]\\
A\ar[r]&X}\in G_P(A)\ ,$$
where $\phi(a,u):=s(a)\bar u$ and $\bar u$ denotes the image of $u\in U$ under $U\to PU$. We consider this object as a morphism $A\to G_P$.
We leave it as an exercise to verify that this map is an atlas.

\subsubsection{}\label{gerbw}

A morphism of stacks $G\rightarrow X$ with $X$ a space is a topological gerbe with band $U(1)$ if there exists an atlas
$A\rightarrow X$, a lift
\begin{equation}\label{sswa}\xymatrix{&G\ar[d]\\A\ar[r]\ar@{.>}[ur]&X}\end{equation}
to an atlas of $G$ such that
$$\xymatrix{U(1)\ar[d]&&\\ A\times_GA\ar@{=>}[r]\ar[d]& A\ar@{=}[d]&\\A\times_XA\ar@{=>}[r]&A\ar[r]&X}$$
is a $U(1)$-extension of topological groupoids.
 In particular, the bundle gerbes considered in \ref{rrtp} give rise to topological gerbes with band $U(1)$.
For equivalent definitions see  \cite{math.AG/0503247}, \cite{heinloth}.
The definition of a gerbe in 
\cite{MR1197353} is slightly more general since the existence of an atlas is not required.

\subsubsection{}

The lifting gerbe $G_P\to X$ constructed in \ref{lg} is a topological gerbe with band $U(1)$. In fact, the construction \ref{atla} produces the  lift (\ref{sswa}).

\subsubsection{}

In the definitions above the Grothendieck site $\Top$ can be replaced by the Grothendieck site of smooth manifolds  $\Mf^\infty$. In this site the covering families are again coverings by families of open submanifolds.

Stacks on $\Mf^\infty$
are called stacks in smooth manifolds.
If $G$ is a stack in smooth manifolds, then an atlas $A\to G$ is a map of stacks which is representable and smooth, i.e. for any map
$T\to G$ from a smooth manifold $T$ to $G$ the induced map
$T\times_GA\to A$ is a submersion of manifolds.
A stack in smooth manifolds which admits an atlas will then be called smooth.

\subsubsection{}
Let $Y\rightarrow X$ be a map of manifolds. It is representable as a map between stacks in smooth manifolds if for any map
$Z\rightarrow X$ the fibre product $Z\times_XY$ exists as a manifold. 
Submersions between manifolds are representable maps.\footnote{ We do not know the converse, i.e. whether  a representable map between manifolds is necessarily a submersion.}

\subsubsection{}

We come to the conclusion that a basic object classified by 
$\lambda_\Z\in H^3(X;\Z)$ is the equivalence  class of a smooth  gerbe $f:G\rightarrow X$ with band
$U(1)$. Instead of going the way through some version of parametrized stable homotopy theory it now seems natural to define a real cohomology twisted by $G$ directly using a suitable  sheaf theory on stacks. A natural candidate would be something like
$H^*(X;G):=H^*(G;\R):=H^*(G;i^\sharp \R_G)$, where
$i^\sharp \R_G$ is the sheafification of the constant presheaf with value $\R$, and $H^*(\dots,i^\sharp \R_G)$ is defined using the derived global sections, or the derived $p_*$, where $p:G\rightarrow *$ is the projection to a point. In fact, if $G$ would be a  manifold, then the sheaf theoretic $H^*(G,i^\sharp\R_G)$ would be isomorphic to the de Rham cohomology of the manifold $G$, and therefore to  the topologist's $H^*(G;\R)$. 

To proceed in the case of stacks we must clarify what we mean by a sheaf on $G$, and how we define $p_*$. The construction of $H^*(G;\R)$ will be finalized in Definition \ref{cogo}. In order to define  sheaves and presheaves on $G$ we associate in \ref{thesite} to $G$ a Grothendieck  site $\bG$. The notions of presheaves and sheaves on a site are the standard ones.

\subsubsection{}

To define cohomology for stacks one can use different sites.
The choices in \cite{math.DG/0605694}
and \cite{heinloth} differ from
our choice, but we indicate that the resulting
cohomologies can be compared and are isomorphic (\ref{site-comp}).
One of our main aims is to study the functorial
properties of the derived categories of sheaves
attached to the sites $\bG$, the functoriality is used here and in subsequent
work, in particular in \cite{bssm}, where we use functoriality to obtain a
periodization with good properties of ordinary cohomology on stacks.

\subsubsection{}

So, if $f:G\rightarrow X$ is a morphism of stacks, then we are interested in functors $f_*,f^*$. Such operations are usually obtaind from some induced morphisms of sites $f^\sharp:\bX\rightarrow \bG$. In fact, this works well for representable morphisms. But in the case of a gerbe
$f:G\rightarrow X$ neither $f$ nor $p:G\rightarrow *$ are representable.
We will define $f_*$ and $p_*$ in an ad-hoc way. The same problem with a similar 
solution also occurs in algebro-geometric set-ups, see e.g. \cite{MR1771927}.
Because of this ad-hoc definitions we must redevelop some of the basic material of sheaf theory in order to check that the expected properties hold in the present set-up. 
For details we refer to the introduction \ref{introsch} to sheaf theory part of the present paper.
 
\subsubsection{}

After the development of elements of sheaf theory on smooth stacks we can define
$$H^*(X;G):=H^*(G;\R):=H^*(\ev\circ Rp_*\circ Ri(i^\sharp \R_G))\ ,$$
where $i:\Sh\bG\rightarrow \Pr\bG$ is the embedding of sheaves into presheaves, the sheafification functor $i^\sharp:\Pr\bG\to \Sh\bG$ is the left adjoint of $i$, and the exact functor 
$$\ev:{\Pr}_{\Ab}\Site(*)\rightarrow \Ab$$ evaluates a presheaf of abelian groups on the object $(*\rightarrow *)\in \Site(*)$. This last evaluation is necessary since our site is the big site of $*$ consisting of all smooth manifolds.  
As the notation suggests we view this as the cohomology of $X$ twisted by the gerbe $G$.

\subsubsection{}
This definition is natural in $G$. If $u:G^\prime\rightarrow G$ is a smooth map of stacks, then by Lemma \ref{pulb} we have a functorial 
map $$u^*:H^*(G;\R)\rightarrow H^*(G^\prime;\R)$$
since there is a canonical isomorphism $u^*i^\sharp \R_G\cong i^\sharp \R_{G^\prime}$. In particular, $H^*(X;G)$
carries the action of the automorphisms of the gerbe $G\rightarrow X$.
One can define the map $u^*$ without the assumption that $u$ is smooth, but then the argument is more complicated, see \cite{bss1}.

\subsubsection{}

The natural question is now how the  $\lambda$-twisted de Rham cohomology $H^*(X;\lambda)$ and $H^*(X;G)$ are related. The main step in this relation is provided by Theorem \ref{main}.
Using this result in the isomorphism $\stackrel{!}{\cong}$ and the projection $q:X\rightarrow *$ we can write
\begin{eqnarray*}
H^*(X;G)&=& H^*(\ev\circ Rp_*\circ Ri(i^\sharp R_G))\\
&\cong&H^*(\ev\circ R(q\circ f)_*\circ Ri(i^\sharp R_G))\\ 
&\stackrel{(**)}{\cong}& H^*(\ev\circ Rq_*\circ Rf_*\circ Ri(i^\sharp \R_G))\\
&\stackrel{(*)}{\cong} & H^*(\ev\circ Rq_*\circ Ri\circ i^\sharp\circ Rf_*\circ Ri(i^\sharp \R_G))\\
&\stackrel{!}{\cong}& 
H^*(\ev\circ Rq_*\circ Ri(\Omega^\cdot[[z]]_\lambda))\\
&\stackrel{(****)}{\cong}&H^*(\ev\circ R(q_*\circ i)(\Omega^\cdot[[z]]_\lambda))\\
&\stackrel{(***)}{\cong}&H^*(\ev\circ q_*\circ i(\Omega^\cdot[[z]]_\lambda))\\
&= &H^*(\Omega^\cdot[[z]]_\lambda(X))\ .\end{eqnarray*}
In order to justify the isomorphism $(*)$ we use Lemma \ref{shpre} which says that
$f_*$ preserves sheaves. The isomorphism $(**)$ follows from Lemma \ref{rcompo} since $f$ is smooth. For $(****)$ we use Lemma \ref{ddddw}. Finally,
$(***)$ follows from Lemma \ref{klom} and the fact that $\Omega^\cdot[[z]]_\lambda$ is a complex of flabby sheaves (see \ref{flabbydef}).

 Note that the isomorphism $\stackrel{!}{\cong}$ depends on additional choices.

\subsubsection{}

It remains to relate the cohomology of the complex
$(\Omega^\cdot[[z]]_\lambda(X),d_\lambda)$ (see \ref{oldes}) with $H^*(X;\lambda)$.
Let $$\Omega^\cdot[[z]]_\lambda(X)^{p}\subset \Omega^\cdot[[z]]_\lambda(X)$$ be the subset of polynomials
$\sum_{2n+k=p} z^n \omega^k$ with $\omega^k\in \Omega^k(X)$. Then we have
$d_\lambda:\Omega^\cdot[[z]]_\lambda(X)^{p}\rightarrow \Omega^\cdot[[z]]_\lambda(X)^{p+1}$.
For $p> 0$ we construct morphisms $\psi_{p}$ such that the following diagram commutes
$$\xymatrix{\dots\ar[r]&\Omega^{odd}(X)\ar[d]^{\psi_{2p-1}}\ar[r]^{d_\lambda}&\Omega^{even}(X)\ar[d]^{\psi_{2p}}\ar[r]^{d_\lambda}&\Omega^{odd}(X)\ar[d]^{\psi_{2p+1}}\ar[r]^{d_\lambda}&\dots\\
\dots\ar[r]&\Omega^\cdot[[z]]_\lambda(X)^{2p-1}\ar[r]^{d_\lambda}&\Omega^\cdot[[z]]_\lambda(X)^{2p}\ar[r]^{d_\lambda}&\Omega^\cdot[[z]]_\lambda(X)^{2p+1}\ar[r]^{d_\lambda}&\dots}$$
In fact, for $e=0,1$ and 
$\omega=\sum_{i=0}^\infty \omega^{e+2i}$ we define 
$$\psi_{2p+e}(\omega):=\sum_{i=0}^{[\frac{p}{2}]} \frac{z^{p-i} \omega^{e+2i}}{(p-i)!}\ .$$
If $p>\dim(X)$, then $\psi_p$ is an isomorphism.
Therefore for large $p$ the isomorphisms $\psi_p$ induce embeddings
$H^*(X;\lambda)\hookrightarrow H^*(X;G)$.
In this way $H^*(X;G)$ is a replacement of $H^*(X;\lambda)$ with good functorial properties.

\subsubsection{}

The definition of real cohomology of $X$ twisted by a gerbe as
$$H^*(X;G):=H^*(G;\R)$$
has a couple of additional interesting  features.
\begin{enumerate}
\item
First of all note that $\R$ is a commutative ring. Therefore
$H^*(X;G)$ has naturally the structure of a graded commutative ring.
In the old picture this structure seems to be partially reflected by the product
$$H^*(X;a\lambda)\otimes H^*(X;b\lambda)\rightarrow H^*(X;(a+b)\lambda)\ .$$
\item
One can replace $\R$ by any other abelian group.
In particular, one can define integral twisted cohomology by
$H^*(X;G;\Z):=H(G;\Z)$.
This definition of an integral twisted cohomology proposes a solution
to the question raised in the remark made in \cite[Sec. 6]{math.KT/0510674}.
Using the maps $\psi_p$ introduced above we can identify the image of
$H^*(X;G;\Z)\rightarrow H^*(X;G)$ as a lattice in $H^*(X;\lambda)$.
The result depends on the choice of $p$, and in view of the denominators in the formula for $\psi_p$ the position of lattice is not very obvious.
\item  
In the proof of Theorem \ref{main} we construct a de Rham model for the cohomology of $H^*(G;\R)$. Let $\Omega_G^{< p}$ be the sheaf of  de Rham complexes truncated at $p-1$ and form the sheaf of Deligne complexes $\cH(p-1)_G:=(i^\sharp \Z_{G}\rightarrow \Omega_G^{<p})$, where 
$i^\sharp \Z_{G}$ sits in degree $-1$.
We can then define the real Deligne cohomology $\hat H^p(G;\Z)$ of $G$ as the $(p-1)$-st hypercohomology
of the complex $\cH_G(p-1)$ (see \cite{MR1197353} for a definition of Deligne cohomology for manifolds in a similar fashion).  
\end{enumerate}

\subsection{Sheaf theory for smooth stacks}\label{introsch}

\subsubsection{}

This subsection is the introduction to the sheaf theoretic part of the paper.
We consider a smooth stack $X$. In order to define the notion of a sheaf on $X$ we associate to $X$ a Grothendieck site $\bX$. In this
paper we adopt the convention of \cite{MR1317816}
that a site consists of a category $\bX$ and the choice of covering families $\cov_\bX(U)$ for the objects $U\in \bX$.
Presheaves on $\bX$ are just contravariant
set-valued functors on $\bX$. A sheaf on $\bX$ is a presheaf which satisfies a descent condition with respect to the covering families.

\subsubsection{}

We define the category $\bX$  as a full subcategory of the category of manifolds $U$ over $X$ such that the structure map $U\to X$ is smooth. The covering families of $U\to X$ are families of submersions over $X$ whose union maps surjectively to $U$.
Observe that the category of smooth manifolds can be considered as a site with the above mentioned choice of covering families. By the Yoneda embedding it maps to the two-category of smooth stacks.
In Subsection \ref{se1} we consider this abstract situation. We consider a site $\cS$, a two-category $\cC$ and a functor $z:\cS\to \cC$. Furthermore we consider a subcatgeory
$r\cC$ which plays the role of the subcategory of stacks with smooth representable morphisms. In this situation we associate to each object $X\in \cC$ the site $\bX$ (see \ref{dd1}) as the full subcategory of $(z(U)\to X)\in \cS/X$ such that the structure map belongs to $r\cC$. The covering families are  induced from $\cS$ (see \ref{dd2}).

\subsubsection{}

The central topic of Subsection \ref{se1} is the adjoint pair (\ref{apari1}) of functors
$$f^*:\Pr\bX\Leftrightarrow \Pr \bG:f_*$$
between presheaf categories associated to a morphism $f:G\to X$. Since in general $f$ does not induce a morphism between the  sites $\bG$ and $\bX$ we define these functors in an ad-hoc manner
(see \ref{dd4} and \ref{dd3}). For two composeable morphisms $f,g$ we relate $(g\circ f)^*, (g\circ f)_*$ with $g_*\circ f_*, f^*\circ g^*$ in Lemma \ref{trtz}.  

\subsubsection{}

In Subsection \ref{thesite} we specialize to smooth stacks.
If the morphism $f:G\to H$  between smooth stacks is smooth or representable, then
it gives rise to a morphism of sites $f_\sharp$ or $f^\sharp$, respectively (\ref{gg} and \ref{gg221}). We verify that our ad-hoc definitions of $f^*$ or
$f_*$ , respectively, coincide with the standard functors induced from the morphism
of sites $f^\sharp$ or $f_\sharp$ (see Lemmas \ref{poi451}, \ref{pi}, \ref{zwei11}).

\subsubsection{}

Most of the statements which we formulate for the sheaf theory on stacks are well-known in
the usual sheaf theory on sites and for functors associated to morphisms of sites.
But for the sheaf theory on stacks we must be very careful about which of these standard facts remain true in general. For other statements  we must know under which additional assumptions they carry over to stacks. 

\subsubsection{}

An important point is the observation that for every morphism between smooth stacks the functor $f_*$ preserves sheaves (Lemma \ref{shpre}). In the  Lemmas \ref{pulb} and \ref{pulb1} we study the compatibility of the pull back with the push forward in cartesian squares. In Lemma \ref{comis} we study under which additional assumptions we have
relations like $(g\circ f)^*\cong f^*\circ g^*$.

\subsubsection{}

In order to define the cohomology of a gerbe we must descend the functors $f_*$ and $f^*$ to the derived categories of presheaves and sheaves of abelian groups. This question is studied  in Subsection \ref{derouzt}. Here the exactness properties of the functors studied in the preceding subsections play an important role. Most of the statements in this subsection are standard for the usual sheaf theory and functors associated to a morphism of sites. Here we study carefully under which additional conditions they remain true for stacks.

\subsubsection{}

The main result (Lemma \ref{schluss}) of Subsection \ref{toola} is that the derived functor  $Rf_*$  for a map $G\to X$ of smooth stack can be calculated using a simplicial approximation of  $G\to X$. In particular, if $X$ is a manifold, then the calculation of $Rf_*$ can  be reduced to ordinary sheaf theory on manifolds.
We use this simplicial model in the proof of our main theorem, for the explicit calculation of the cohomology of the stack  $[*/S^1]$ in Lemma \ref{bas13}, but also to verify that
pull-back and push-forward commute on the level of derived functors for certain cartesian diagrams in Lemma \ref{derpulb}.

\subsubsection{}

The covering families of  the small site  $(U)$ of a manifold are coverings by open subsets.
Thus the sheaf theory for $(U)$ is the ordinary one.  If $(U\to X)\in \bX$, then a presheaf on $X$ induces a presheaf  on $(U)$. In the present paper the sheaf theory on $(U)$ is considered to be well-understood. The main goal of Subsection \ref{cckkll11} is to compare the sheafification functors on $\bX$ and $(U)$ (see Lemma \ref{gloloc}). 
This result is very useful in explicit calculations since it says that certain questions can be studied  for each $(U\to X)\in \bX$ separately and with respect to the small site $(U)$.
This sort of reasoning will be applied in the proof that the de Rham complex of a stack is a flabby resolution of the constant sheaf with value $\R$, where we use that this fact is well-known on each manifold equipped with the site $(U)$.
It is also used in the proof of  Lemma \ref{somecomu} which says that for a smooth map between smooth stacks  the pull back commutes with the sheafification functor.

\section{Sheaf theory for smooth stacks}\label{se1}

\subsection{Over sites}\label{overs}

\subsubsection{}\label{situ1}

The goal of the present subsection is to develop some elements of sheaf theory
in the following situation. Let $\cS$ be a site (see \cite[Chapter I, 1.2.1]{MR1317816} for a definition), $\cC$ a two-category with  invertible two-morphisms, and $z:\cS\rightarrow \cC$ a functor (we consider $\cS$ as a two-category with only identity two-isomorphisms). Finally let  $r\cC$ be a subcategory
of the category underlying $\cC$
which we call the category of admissible  morphisms.

To each object $G\in \cC$ we will associate a site $\bG$ (sometimes we will write $\Site(G):=\bG$) and the categories of presheaves $\Pr\bG$ and sheaves $\Sh\bG$ of sets on this site. For a morphism
$f\in \cC(G,H)$ we will define an adjoint pair of functors
$$f^*:\Pr\bH\Leftrightarrow \Pr\bG:f_*\ .$$
In general these functors are not induced by a morphism of sites.

\subsubsection{}

Let $G\in \cC$. We define the underlying category of $\bG$.
\begin{ddd}\label{dd1}
The objects of $\bG$ are pairs $(U,\phi)$, where $U\in \cS$ and
$\phi\in r\cC(z(U),G)$. A morphism  $(U,\phi)\rightarrow (U^\prime,\phi^\prime)$ is given by a pair $(h,\sigma)$, where $h\in \cS(U,U^\prime)$ and  $\sigma$ is a two-isomorphism 
$$\xymatrix{z(U)\ar[dr]^{\phi}\ar[rr]^{z(h)}&\ar@{=>}[d]^\sigma&z(U^\prime)\ar[dl]^{\phi^\prime}\\&G&}$$
The composition in $\bG$ is defined in the obvious way.
\end{ddd}
Sometimes we will abbreviate the notation and write $U$ or $(z(U)\rightarrow G)$ for $(U,\phi)$.

\subsubsection{}\label{t5tt5t}
Next we define the coverings of an object $(U,\phi)$ of $\bG$. 
\begin{ddd}\label{dd2}
A covering of $(U,\phi)$ is a  collection
of morphisms
$$((U_i,\phi_i)\stackrel{(h_i,\sigma_i)}{\rightarrow} (U,\phi))_{i\in I}$$
such that $(U_i\stackrel{h_i}{\rightarrow} U)_{i\in I}$ is a covering of $U$ in 
$\cS$.
 \end{ddd}

In fact it is easy to verify the axioms listed in the definition 
\cite[1.2.1]{MR1317816}. The only non-obvious part asserts that given a covering
$((U_i,\phi_i)\rightarrow (U,\phi))_{i\in I}$
and a morphism $(V,\psi)\rightarrow (U,\phi)$, then the fibre products
$(V_i,\psi_i):=(V,\psi)\times_{(U,\phi)} (U_i,\phi_i)$ exist in $\bG$ and
$((V_i,\psi_i)\rightarrow (V,\psi))_{i\in I}$ is a covering of $(V,\psi)$.
By a little diagram chase one verifies that
$(V,\psi)\times_{(U,\phi)} (U_i,\phi_i)\cong (V\times_UU_i,\phi\circ z(\kappa))$, where
$V\times_UU_i$ is the fibre product in $\cS$ and $\kappa: V\times_UU_i\rightarrow U$ is the natural map.

\subsubsection{}

Let $f\colon G\rightarrow H$ be a morphism in $\cC$. Then we can define
the functor $f^*:\Pr\bH\rightarrow \Pr\bG$ as follows.
Given $(z(V)\to G)\in \bG$ we consider the category $V/\bH$ 
(recall that $V$ abbreviates $(z(V)\to G)$)
of diagrams
$$\xymatrix{V\ar[d]&&z(V)\ar[r]\ar[d]&G\ar[d]^f\\U&&z(U)\ar@{=>}[ur]\ar[r]&H}\ .$$
A morphism in this category is given by a morphism $(z(U^\prime)\to H)\rightarrow (z(U)\to H)$ in $\bH$  fitting into  
$$\xymatrix{&U&z(V)\ar[dr]\ar[ddr]\ar[rr]&&G\ar[d]^f\\V\ar[ur]\ar[dr]&&&z(U)\ar@{=>}[dr]\ar@{=>}[ur]\ar[r]&H\ar@{=}[d]\\&U^\prime\ar[uu]&&z(U^\prime)\ar[r]\ar[u]&H}\ .$$
Let $F\in \Pr\bH$.
\begin{ddd}\label{dd4}
We define 
$$f^*F(V):=\colim_{V/\bH} F(U)\ .$$
\end{ddd}
A morphism $V^\prime\rightarrow V$ in $\bG$ induces naturally a functor
$V/\bH\rightarrow V^\prime/\bH$. The relevant diagram is
$$\xymatrix{V^\prime\ar[d]&&z(V^\prime)\ar[d]\ar[r]&G\ar@{=}[d]\\V\ar[d]&&z(V)\ar@{=>}[ur]\ar[r]\ar[d]&G\ar[d]^f\\U&&z(U)\ar@{=>}[ur]\ar[r]&H}\ .$$ We therefore get a map
$f^*F(V)\rightarrow f^*F(V^\prime)$, and this makes
$f^*F$ a presheaf on $\bG$.

\subsubsection{}\label{pushdef}

Let $f:G\rightarrow H$ again be a morphism in $\cC$. We define a functor
$f_*:\Pr\bG\rightarrow \Pr\bH$ as follows. We consider
$(z(U)\rightarrow G)\in \bH$. Then we consider the category of diagrams
$\bG/U$ of diagrams
$$\xymatrix{V\ar[d]&&z(V)\ar[d]\ar[r]&G\ar[d]^f\\U&&z(U)\ar@{=>}[ru]\ar[r]&H}\ .$$
A morphism of such diagrams is given by a morphism
$V^\prime\rightarrow V$ in $\bG$ which fits into
$$\xymatrix{V^\prime\ar[d]&&z(V^\prime)\ar[d]\ar[r]&G\ar@{=}[d]\\V\ar[d]&&z(V)\ar@{=>}[ur]\ar[r]\ar[d]&G\ar[d]^f\\U&&z(U)\ar@{=>}[ur]\ar[r]&H}\ .$$
\begin{ddd}\label{dd3}
We define
$$f_*F(U):=\lim_{\bG/U} F(V)\ .$$
\end{ddd}
A morphism $U\rightarrow U^\prime$ in $\bH$ induces naturally a functor
$\bG/U\rightarrow \bG/U^\prime$. The relevant diagram is
$$\xymatrix{V\ar[d]&&z(V)\ar[d]\ar[r]&G\ar@{=}[d]\\U\ar[d]&&z(U)\ar@{=>}[ur]\ar[r]\ar[d]&G\ar[d]^f\\U^\prime&&z(U^\prime)\ar@{=>}[ur]\ar[r]&H}\ .$$
We therefore get a map $f_*F(U^\prime)\rightarrow f_*F(U)$, and this makes
$f_*F$ a presheaf on $\bH$.
 
\subsubsection{}

Let $f\in \cC(G,H)$ as before.
\begin{lem}\label{apari1}
The functors $f_*$ and $f^*$ naturally form an adjoint pair
$$f^*\:\::\:\:\Pr\bH\Leftrightarrow \Pr\bG\:\::\:\:f_*\ .$$
\end{lem}
\proof
We give the unit and the counit.
Let $(z(W)\rightarrow G)\in \bG$. Then
$$f^*f_*F(W)=\colim_{U} \lim_V F(V)\ ,$$
where the colimit-limit is taken over a category of diagrams
$$\xymatrix{V\ar[d]&&z(V)\ar[d]\ar[r]&G\ar[d]^f\\U&&z(U)\ar[r]&H\\W\ar[u]&&z(W)\ar[u]\ar[r]&G\ar[u]_f}$$
(we leave out the two-isomorphisms).
The counit is a natural transformation 
$$f^*f_*F(W)\rightarrow F(W)\ .$$
It is given by the universal property of the colimit and the collection of maps which associates to $U$ the canonical map $\lim_V F(V)\to F(W)$.

Furthermore, let $(z(U)\rightarrow H)\in \bH$.
Then
$$f_*f^*F(U)=\lim_{V}\colim_{W} F(W)\ ,$$
where the limit-colimit is taken over a category of diagrams
$$\xymatrix{U&&z(U)\ar[r]&H\\V\ar[u]\ar[d]&&z(V)\ar[u]\ar[d]\ar[r]&G\ar[u]_f\ar[d]^f\\W&&z(W)\ar[r]&H}$$
(we leave out the two-isomorphisms).
The unit is a natural transformation
$$F(U)\rightarrow f_*f^*F(U)\ .$$ 
It is given by the universal property of the limit and the collection of maps
which associates to $V$ the natural map
$F(U)\to \colim_W F(W)$.

We leave it to the interested reader to perform the remaining checks.
\hB

\subsubsection{}

Let us consider a pair of  composable maps in $\cC$
$$G\stackrel{f}{\longrightarrow} H\stackrel{g}{\longrightarrow}L\ .$$
\begin{lem}\label{trtz}
We have natural transformations of functors
$$(g\circ f)_*\to g_*\circ f_*\ ,\quad f^*\circ g^*\to (g\circ f)^*\ .$$
\end{lem}
\proof
We discuss the transformation $f^*\circ g^*\to (g\circ f)^*$.
Let $F\in \Pr\bL$ and $(z(W)\to G)\in \bG$. Inserting the definitions we have 
$$f^*\circ g^*(F)(W)\cong \colim_{\cA} F(U)\ ,$$
where $\cA$ is the  category of diagrams
$$\xymatrix{z(W)\ar[r]\ar[d]&G\ar[d]^f\\
z(V)\ar[r]\ar[d]\ar@{=>}[ur]&H\ar[d]^g\\
z(U)\ar[r]\ar@{=>}[ur]&L}$$ with $(V\to H)\in \bH$ and $(U\to L)\in \bL$.
The vertical composition
provides a functor $\cA\to W/\bL$, where $W/\bL$ is the category  of diagrams of the form
$$\xymatrix{z(W)\ar[d]\ar[r]&G\ar[d]^{g\circ f}\\
z(U)\ar[r]\ar@{=>}[ur]&L}\ .$$
We get an induced map of colimits
$$f^*\circ g^*(F)(W)\to (g\circ f)^*F(W)=\colim_{W/\bL}F(U)\ .$$
The other transformation $(g\circ f)_*\to g_*\circ f_*$ is obtained in a similar manner or, equivalently, by adjointness.
\hB

In general, we can not expect that these transformations are isomorphisms. But under additional assumptions they are, see \ref{comis}.

\subsection{The site of a smooth stack}\label{thesite}

\subsubsection{}

We consider the site $\Mf^\infty$ of smooth manifolds\footnote{In order to avoid set-theoretic problems one must require that a site is a small category. In the present paper we will ignore this problem. It can be resolved by either working with universes or replacing $\Mf^\infty$  by an equivalent small category (see e.g. \cite{math.DG/0306176}).} and open covering families.
Its underlying category is the category of smooth manifolds and smooth maps. A collection of smooth maps 
$(U_i\rightarrow U)_{i\in I}$ is a covering if and only if this family is isomorphic to
the collection of inclusions of the open subsets of $U$ given by an open covering of $U$. 

We use the site $\Mf^\infty$ in order to define stacks in smooth manifolds.
We refer to \cite{heinloth}, \cite{math.DG/0306176}, \cite{math.AG/0503247} for the language of stacks.

\subsubsection{}

We will also consider the site $\cS$ on smooth manifolds.
In this site a family $(U_i\rightarrow U)_{i\in I}$ of smooth maps is a covering
if the maps $U_i\rightarrow U$ are submersions and $\sqcup_{i\in I} U_i\rightarrow U$ is surjective. We will use this site in order to define the site of a stack according to \ref{overs}.
In fact the descent conditions for $\Mf^\infty$ and $\cS$
are the same, and it is only a matter of taste
that we use the notion site in this way.

\subsubsection{}\label{site33}
In this paragraph we recall the main notions of the theory of smooth stacks.
\begin{enumerate}
\item
A morphism of stacks $G\rightarrow H$ is called representable, if
for each manifold $U$ and map $U\rightarrow H$ the fibre product $U\times_HG$ is equivalent to a manifold. 
A composition of representable maps is representable.
\item
A representable morphism $G\rightarrow H$ of stacks is called smooth if for each manifold $U$ and map $U\to H$ the induced map $U\times_HG\to U$ (of manifolds) is a submersion.
\item
A map $U\to G$ from a manifold to a stack is called an atlas if it is representable, smooth and admits local sections. 
\item
A stack in smooth manifolds is called smooth if it admits an atlas \cite[Def. 2.4]{heinloth}. 
\item
A morphism (not necessarily representable) between smooth stacks $f:G\to H$ is called smooth if for an atlas $A\to G$ the composition
$A\to G\to H$ is smooth \cite[Def. 2.10]{heinloth}. 
A composition of smooth maps is smooth.
\item 
A smooth morphism $U\to G$ from a manifold to a smooth stack is representable.
\end{enumerate}

\subsubsection{}\label{ccea}

Let $\cC$ be the two-category of smooth stacks in smooth manifolds.  
We have  a Yoneda embedding $z:\cS\rightarrow \cC$. Note that in general we will omit the Yoneda embedding in the notation and consider $\cS$ as a subcategory of $\cC$. We let $r\cC$ be the subcategory of representable smooth morphisms.

\subsubsection{}

The conventions introduced in \ref{ccea} place us in the situation of \ref{situ1}.
Let $G\in\cC$ be a smooth stack. Then by  $\bG$  we denote the site according
to Definitions \ref{dd1} and \ref{dd2}. Note that this site is derived from
the site $\cS$ on smooth manifolds. We now have the categories of presheaves
$\Pr\bG$ and sheaves $\Sh\bG$ on the stack $G$. We compare the present
definition to the one of \cite[Sec. 4]{heinloth} in Subsection
\ref{site-comp}. % execpt that we use the big site $\cS$ on the category of smooth manifolds instead of the small site $\Mf^\infty$.

\subsubsection{}\label{gg}

Let $f:G\rightarrow H$ be a representable morphism of smooth stacks in smooth manifolds.
Then it induces a morphism of sites $f^\sharp:\bH\rightarrow \bG$
 by the rule $f^\sharp(U\rightarrow H):=U\times_H G\rightarrow G$ (it is easy to check the axioms listed in \cite[1.2.2]{MR1317816}). 

\subsubsection{}\label{gg221}

If $f:G\to H$ is a smooth morphism of smooth stacks in smooth manifolds, then we can define another morphism of sites
$f_\sharp:\bG\to \bH$ by $f_\sharp(V\to G):=(V\to G\stackrel{f}{\to} H)$.

\subsubsection{}

We call a functor left exact if it preserves arbitrary limits. If it preserves arbitrary colimits, then we call it right exact. A functor is said to be exact if it is right and left exact.

Recall that a functor which is a left adjoint is right exact. Similarly, a right adjoint is left exact.

\subsubsection{}

A morphism of sites
$q:\bH\rightarrow \bG$ induces an adjoint pair
$$q_*\:\::\:\:\Pr\bH\Leftrightarrow \Pr\bG\:\::\:\:q^*\ .$$
(see \cite[2.3]{MR1317816}).
In the following we compare these maps with the ad-hoc  definitions \ref{dd4} and \ref{dd3}  and discuss some special properties.

\subsubsection{}

\begin{lem}\label{poi451}
If $f:G\to H$ is a smooth morphism between smooth stacks, then
we have $f^*\cong (f_\sharp)^*$. In particular, then $f^*$ is exact and preserves sheaves.
\end{lem}
\proof
Let $(V\to G)\in \bG$.
According to the definition \cite[2.3]{MR1317816} we have
$$(f_\sharp)^*(F)(V\to G):=F(f_\sharp(V\to G))=F(V\to G\to H)\ .$$ 
If $(V\to G)\in \bG$, then the category
$V/\bH$ has an initial object
$$\xymatrix{V\ar[r]\ar@{=}[d]&G\ar[d]^f\\V\ar@{=>}[ur]^{\id}\ar[r]&H}\ .$$
Therefore
\begin{equation}\label{ttgg221}(f^*F)(V\to G)\cong F(V\to G\to H)\ .\end{equation}
This implies that $f^*\cong (f_\sharp)^*$.

It is well-known \cite[3.6]{MR1317816} that 
the contravariant functor (in our case $(f_\sharp)^*$) associated to a morphism of sites preserves sheaves. Therefore $f^*$ preserves sheaves.

The limit of a diagram of presheaves is defined objectwise. By  (\ref{ttgg221}) the functor $f^*$ commutes with limits. As a left adjoint (by Lemma \ref{apari1}) it also commutes with colimits.
\hB

\subsubsection{}

Let $f:G\rightarrow H$ be a representable and smooth morphism of smooth stacks.
 
\begin{lem}\label{pi}
We have an isomorphism of functors $(f^\sharp)_*\cong f^*$.
\end{lem}
\proof
Let $F\in \Pr\bH$.
For $(V\rightarrow G)\in \bG$ we have the category
$V/f^\sharp$ of pairs $((U\rightarrow H)\in \bH,(V\rightarrow f^\sharp(U))\in \Mor(\bG))$.
It has a natural evaluation $\ev_V:V/f^\sharp\rightarrow \bH$ which maps
$((U\rightarrow H),(V\rightarrow f^\sharp(U)))$ to $(U\to H)$.
By definition (see \cite[Proof of 2.3.1]{MR1317816})
$$(f^\sharp)_*(F)(V)=\colim_{V/f^\sharp} F\circ \ev_V\ .$$

Now we observe that $V/f^\sharp$ can be identified with the category of diagrams
$$\xymatrix{V\ar[d]\ar[r]&U\ar[d]\\
G\ar[r]^f\ar@{=>}[ur]&H}\ .$$
Since $f$ is smooth we see that
$(V\to G\stackrel{f}{\to} H)\in \bH$ and 
$$\xymatrix{V\ar[d]\ar[r]^{\id_V}&V\ar[d]\\
G\ar[r]^f\ar@{=>}[ur]^{\id}&H}$$
is the initial element of $V/f^\sharp$. We conclude that
\begin{equation}\label{vali}(f^\sharp)_*(F)(V\to G)\cong F(V\to G\to H)\ .\end{equation}
The equality $f^*\cong (f^\sharp)_*$ now follows from  (\ref{ttgg221}).
\hB

One can not expect that $f^*$ is left exact for a general map $f:G\rightarrow H$. In fact this problem occurs in the corresponding definition in \cite{MR1771927} of the pull-back for the lisse-etale site of an algebraic stack. For more details and a solution see \cite{olsson}.

% \subsubsection{}

% \begin{lem}
% If $f:G\rightarrow H$ is a representable and smooth morphism of smooth stacks, then we have a natural isomorphism
% $f^*\cong (f^\sharp)_*$. 
% \end{lem}
% \proof
% If $f$ is representable and $(V\to G)\in \bG$, then the category
% $V/\bH$ has an initial object
% $$\xymatrix{V\ar[r]\ar@{=}[d]&G\ar[d]^f\\V\ar@{=>}[ur]^{\id}\ar[r]&H}\ .$$
% Therefore for $F\in \Pr\bH$ we have  $f^*F(V\to G)=F(V\to G\to H)$.
% Comparison with (\ref{vali}) gives the assertion.
% \hB 

\subsubsection{}

\begin{lem}\label{zwei11}
If $f:G\rightarrow H$ is a representable morphism of smooth stacks, then
$f_*=(f^{\sharp})^*:\Pr \bG\rightarrow \Pr\bH$. The functor $f_*$ is exact.
\end{lem} 
\proof
Let $(U\to H)\in \bH$. Then
$f^\sharp(U\to H)=(U\times_HG\to G)$ is the final object in $\bG/U$. Therefore
\begin{equation}\label{objw2}f_*F(U)\cong F(U\times_HG)\cong (f^\sharp)^*F(U)\ .\end{equation} 
Since $(f^{\sharp})^*$ is a right adjoint it commutes with limits.
Since colimits of presheaves are defined objectwise it follows from the formula
(\ref{objw2}) that $f_*$ also commutes with colimits.
\hB

\subsubsection{}

Let now $f:G\rightarrow H$ be a map of smooth stacks.

\begin{lem}\label{shpre}
The functor $f_*$ preserves sheaves.
\end{lem}
\proof
Let $F\in \Sh\bG$. Consider $(U\rightarrow H)\in \bH$ and let $(U_i\rightarrow U)$ be a covering of $U$. Consider a diagram
\begin{equation}\label{piu}\xymatrix{V\ar[d]\ar[r]&G\ar[d]^f\\U\ar@{=>}[ru]\ar[r]&H}\ .\end{equation}
>From this we obtain a collection of diagrams
$$\xymatrix{V_i:=U_i\times_UV\ar[d]\ar[r]&G\ar[d]^f\\U_i\ar@{=>}[ru]\ar[r]&H}$$functorially in $V$. Observe that $(V_i\rightarrow V)$ is a covering in $\bG$.
We now consider the map of diagrams
$$\xymatrix{f_*F(U)\ar[d]\ar[r]&\prod_i f_*F(U_i)\ar[d]\ar@{=>}[r]&\prod_{i,j}f_*F(U_i\times_UU_j)\ar[d]\\F(V)\ar[r]&\prod_i F(V_i)\ar@{=>}[r]&\prod_{i,j}F(V_i\times_VV_j)
}\ .$$
The vertical maps are given by specialization.
We must show that the upper horizontal line is an equalizer diagram.
The lower horizontal line has this property since $F$ is a sheaf.

We now take the limit over the category of diagrams
(\ref{piu}) und use the fact that a limit preserves equalizer diagrams.
We get the commutative diagram of sets
$$\xymatrix{f_*F(U)\ar@{=}[d]\ar[r]&\prod_i f_*F(U_i)\ar[d]^s\ar@{=>}[r]&\prod_{i,j}f_*F(U_i\times_UU_j)\ar[d]\\f_*F(U)\ar[r]&\lim \prod_i F(V_i)\ar@{=>}[r]&\lim \prod_{i,j}F(V_i\times_VV_j)
}\ .$$
Let us  assume that $s$ is injective. 
Then the fact that the lower horizontal line is an equalizer diagram 
implies by a simple diagram chase
that the upper horizontal line is an equalizer diagram.

We now show that $s$ is injective.
Note that a priori the product of specialization maps
$$s=\prod_i s_i : \prod_i f_*F(U_i)\to \prod_i \lim F(V_i)$$
may  not be injective since the functors
$L_i:\bG/U\ni V\mapsto V_i\in \bG/U_i$ are not necessarily essentially surjective.
But in our situation the maps $s_i$ are injective since each object in $\bG/U_i$ maps into an object in the image of $L_i$.
To see this consider 
a diagram 
$$\xymatrix{W\ar[d]\ar[r]^t& G\ar[d]^f\\U_i\ar@{=>}[ur]\ar[r]&H}\in \bG/U_i\ .$$ Using the composition $W\to U_i\to U$ we can form the diagram
$$\xymatrix{W\ar[d]\ar[r]& G\ar[d]^f\\U\ar@{=>}[ur]\ar[r]&X}\in \bG/U$$
and define the morphism in $\bG/U_i$
$$\xymatrix{U_i\times_U W\ar@/_1cm/[dd]_{\pr_1}\ar[dr]^{t\circ \pr_2}&\\W\ar[u]^j\ar[d]\ar[r]& G\ar[d]^f\\U_i\ar@{=>}[ur]\ar[r]&H}\ ,$$
where $j:W\to U_i\times_UW$ is induced by $W\to U_i$ and $\id_W:W\to W$.

\hB 

\subsubsection{}

Assume that we have a diagram in smooth stacks
\begin{equation}\label{squ567}\xymatrix{G\ar[r]^{u}\ar[d]^f&H\ar[d]^{g}\\
M\ar@{=>}[ur]\ar[r]^v&N}\ ,\end{equation}
 where $u$  and $v$ are smooth.

\begin{lem}\label{pulb}
We have a natural map of functors $\Pr\bH\rightarrow \Pr\bM$
$$v^*\circ g_*\rightarrow f_*\circ u^*$$
which is an isomorphism if
(\ref{squ567}) is cartesian. 
\end{lem}
\proof
We use the description (\ref{ttgg221})  of $v^*$ obtained in the proof of Lemma \ref{poi451}.
Let $(U\to M)\in \bM$. Then we have
$$v^*\circ g_* F(U)\cong \lim F(A)\ ,$$
where the limit is taken over a category $D$ of diagrams
\begin{equation}\label{diq1}\xymatrix{&&&A\ar[dl]\ar[ddd]&\\&G\ar[r]^{u}\ar[d]^f&H\ar@{=>}[ddr]^{\sigma}\ar[d]^{g}&\\
&M\ar@{=>}[drr]\ar@{=>}[ur]\ar[r]^v&N&\\
U\ar[ru]\ar[rrr]^{\id}&&&U\ar[ul]}\end{equation}
(where $A$ varies).
On the other hand
$f_*\circ u^*(F)(U)\cong \lim F(V)$, where the limit is taken over the category $E$ of diagrams
\begin{equation}\label{ee122}
\xymatrix{V\ar[dr]\ar[ddd]\ar[rrr]&&&V\ar[dl]&\\&G\ar@{=>}[ddl]\ar@{=>}[urr]\ar[r]^{u}\ar[d]^f&H\ar[d]^{g}&\\
&M\ar@{=>}[ur]\ar[r]^v&N&\\
U\ar[ru]&&&}
\end{equation}
($V$ varies).
We define a functor
$X:E\to D$ which sends the diagram (\ref{ee122}) to the diagram
$$\xymatrix{&&&V\ar[ddd]\ar[dl]&\\&G\ar[r]^{u}\ar[d]^f&H\ar@{=>}[ddr]\ar[d]^{g}&\\
&M\ar@{=>}[ur]\ar@{=>}[drr]\ar[r]^v&N&\\
U\ar[ru]\ar[rrr]^{\id}&&&U\ar[ul]}\ .$$
We write $F_E$ and $F_D$ for the functor $F$ precomposed with the evaluations
$E\to \bH$ and $D\to \bH$. The identity $F(V)\stackrel{\sim}{\to} F(V)$ induces an isomorphism
$F_D\circ X\stackrel{\sim}{\to} F_E$. 
Therefore we have a natural map of limits
\begin{equation}\label{hin123}
v^*\circ g_* F(U)\rightarrow f_*\circ u^*F(U)\ .
\end{equation}
This gives the required transformation of functors.

If (\ref{squ567}) is cartesian, then we can define a functor
$Y:D\to E$ which maps the diagram (\ref{diq1}) functorially to 
$$\xymatrix{U\times_NA\ar[dr]^{h}\ar[ddd]&&&&\\&G\ar@{=>}[ddl]\ar[r]^{u}\ar[d]^f&H\ar[d]^{g}&\\
&M\ar@{=>}[ur]\ar@{=>}[drr]\ar[r]^v&N&\\
U\ar[ru]\ar[rrr]^{\id}&&&U\ar[ul]}\ ,$$
which employs the map $A\to H\to N$.
The map $h$ is induced by the universal property of the cartesian diagram.
Since $U\to M$ and $A\to H$ are smooth, the map $h$ is smooth, too.
The map $A\to U$ together with the two-isomorphism $\sigma$ gives a map
$A\to U\times_NA$ in $\bH$. This map induces the natural transformation
$F_E\circ Y \to F_D$. It  gives a map of limits
 \begin{equation}\label{her123}
f_*\circ u^*(F)(U)\to v^*\circ g_*F(U)\ .
\end{equation}
One can check that (\ref{her123}) is inverse to (\ref{hin123}).
\hB

\subsubsection{}

Assume again that we have a diagram in smooth stacks
\begin{equation}\label{squ}\xymatrix{G\ar[r]^{u}\ar[d]^f&H\ar[d]^{g}\\
M\ar@{=>}[ur]\ar[r]^v&N}\ .\end{equation}
We now assume that $f$ and $g$ are representable, and that $u,v$ are smooth.
 \begin{lem}\label{pulb1}
We have a natural map of functors $\Pr\bH\rightarrow \Pr\bM$
$$v^*\circ g_*\rightarrow f_*\circ u^*$$
which is an isomorphism if
(\ref{squ}) is cartesian. 
\end{lem}
This is a special case of Lemma \ref{pulb}. But under the additional representablility assumptions on $f$ and $g$ the proof simplifies considerably. 

\proof
Let $F\in \Pr \bH$. For $(U\to M)\in \bM$ we calculate
\begin{eqnarray*}
v^*\circ g_*(F)(U)&\cong& \colim_{V\in U/\bN} \lim_{W\in \bH/V} F(W)\\
&\stackrel{(\ref{ttgg221})}{\cong}&\lim_{W\in \bH/(U\to M\to N)} F(W)\\
&\stackrel{(\ref{objw2})}{\cong}&F(H\times_NU\to H)\ .\end{eqnarray*}
On the other hand
\begin{eqnarray*}f_*\circ u^*(F)(U)&\cong&\lim_{Z\in \bG/U}\colim_{W\in Z/\bH} F(W)\\
&\stackrel{(\ref{objw2})}{\cong}&\colim_{W\in (U\times_MG)/\bH} F(W)\\
&\stackrel{(\ref{ttgg221})}{\cong}&F(U\times_MG\to G\to H)\ .
\end{eqnarray*}
The transformation
$v^*\circ g_*(F)(U)\to f_*\circ u^*(F)(U)$ is now induced from the map $(G\times_MU\to H\times_N U)\in \bH$.

If the diagram is cartesian, then we have
$G\times_MU\cong (H\times_NM)\times_MU\cong H\times_NU$ so that the transformation is an isomorphism.
\hB

\subsubsection{}

Let us consider a pair of  composable maps of smooth stacks
$$G\stackrel{f}{\longrightarrow} H\stackrel{g}{\longrightarrow}L\ .$$
In Lemma \ref{trtz} we have found natural transformations of functors between presheaf categories $$(g\circ f)_*\to g_*\circ f_*\ ,\quad f^*\circ g^*\to (g\circ f)^*\ .$$
\begin{lem}\label{comis}
If $g$ is representable, or if $f$ is smooth, then these transformations are isomorphisms.
\end{lem}
\proof
We consider the transformation $f^*\circ g^*\to (g\circ f)^*$ which appears as a transformation of colimits induced by a functor between indexing categories $\cA\to W/\bL$, where
we use the notation introduced in the proof of Lemma \ref{trtz}.
Under the present additional assumptions on $f$ or $g$ we have a functor $W/\bL\to \cA$
which induces the inverse of the transformation. In the following we describe these functors.

If $g$ is representable, then
each diagram
\begin{equation}\label{rrwq1}\xymatrix{W\ar[d]\ar[r]&G\ar[d]^{g\circ f}\\
U\ar[r]\ar@{=>}[ur]&L}\end{equation}
in $W/\bL$ naturally completes to
$$\xymatrix{W\ar[r]\ar[d]&G\ar[d]^f\\
U\times_LH\ar[r]\ar[d]\ar@{=>}[ur]&H\ar[d]^g\\
U\ar[r]\ar@{=>}[ur]&L}$$ 
in $\cA$.

If $f$ is smooth,
then the diagram (\ref{rrwq1}) can be naturally completed to
$$\xymatrix{W\ar[r]\ar@{=}[d]&G\ar[d]^f\\
W\ar[r]\ar[d]\ar@{=>}[ur]^{\id}&H\ar[d]^g\\
U\ar[r]\ar@{=>}[ur]&L}$$ 
in $\cA$.

It follows from adjointness that under the additional assumptions on $f$ or $g$ the transformation $(g\circ f)_*\to 
g_*\circ f_*$ is an isomorphism, too.
\hB

\subsubsection{}

Let $f:G\rightarrow H$ be a smooth map of smooth  stacks.
The following Lemma is standard, we include
a proof for the sake of completeness.

\begin{lem}\label{lef}
There exists a functor $f_!:\Pr\bG\to \Pr\bH$ so that we get an adjoint pair
$$f_!:\Pr\bG\Leftrightarrow \Pr\bH:f^*\ .$$
\end{lem}
\proof
Let $(V\to G)\in \bG$. Then by (\ref{ttgg221}) we have
 $f^*F(V)\cong F(V\to G\to H)$.

Let $(V\to G)\in \bG$ and $h_{V\to G} \in \Pr\bG$ be the corresponding representable presheaf. Then we have a natural isomorphism
\begin{eqnarray*}
\Hom(h_{V\to  G},f^*F)&\cong&f^*F(V\to G)\\
&\cong&F(V\to G\to H)\\
&\cong&\Hom(h_{V\to G\to H},F)\end{eqnarray*}
which leads us to the definition 
$$f_! h_{V\to G}:=h_{V\to G\to H}\ .$$
If $L\in \Pr\bG$, then we can write
$L\cong \colim_{h_{V\to G}\to L}h_{V\to G}$. Since a left-adjoint must commute with colimits we are forced to set
$$f_!L:=\colim_{h_{V\to G}\to L}h_{V\to G\to H}\ .$$
Then we have indeed
\begin{eqnarray*}
\Hom(L,f^*F)&\cong&\Hom(\colim_{h_{V\to G}\to L}h_{V\to G},f^*F)\\
&\cong&\lim_{h_{V\to G}\to L} \Hom(h_{V\to G},f^*F)\\
&\cong&\lim_{h_{V\to G}\to L} \Hom(h_{V\to G\to H},F)\\
&\cong&\Hom(\colim_{h_{V\to G}\to L}h_{V\to G\to H},F)\\
&\cong&\Hom(f_!L,F)
\end{eqnarray*}
 \hB

\subsection{Presheaves of abelian groups and derived functors}\label{derouzt}

\subsubsection{}

In the previous Subsections \ref{overs} and \ref{thesite} we have developed a theory of set-valued presheaves and sheaves on stacks. We are in particular interested in the abelian categories of presheaves and sheaves of abelian groups and their derived categories. The functors $(f^*,f_*)$ and $(i^\sharp,i)$ preserve abelian group valued objects. In the present subsection we study how these functors descend to the derived categories. Furthermore, we check some functorial properties of these descended functors 
which will be employed in later calculations.

The derived version (Lemma \ref{derpulb}) of the fact that pull-back commutes with push-forward in certain cartesian diagrams (Lemma \ref{pulb}) would fit into the present subsection, but can only be shown after the development of a computational tool in Subsection \ref{toola}.

A similar remark applies to  Lemma \ref{somecomu} saying  that sheafification commutes with
pull-back along smooth maps between smooth stacks. We will show this Lemma in 
Subsection \ref{cckkll11}.

\subsubsection{}

For a site $\bG$ let
$\Pr_{\Ab}\bG$ and $\Sh_{\Ab}\bG$ denote the abelian categogies of presheaves and sheaves of abelian groups on $\bG$. These categories have enough injectives \cite[2.1.1 and 2.1.2]{MR1317816}.
Let $D^+(\Pr_{\Ab}\bG)$ and $D^+(\Sh_{\Ab}\bG)$ denote the lower bounded derived categories of $\Pr_{\Ab}\bG$ and $\Sh_{\Ab}\bG$. 
% Recall that a presheaf of abelian groups on a stack $G$ is a
% functor $F:\cat(\bG)^{op}\rightarrow \Ab$. It is a sheaf, if for all
% $(U,f)\in \cat(\bG)$ und $(U_i\rightarrow U)_{i\in I}\in \cov(\bG)$ the sequence
% $$0\rightarrow F(U)\rightarrow \sqcup_{i\in I}F(U_i)\rightarrow \sqcup_{(i,j)\in I\times I}F(U_i\times_UU_j)$$ is exact.

% We let $\Pr\bG$ and $\Sh\bG$ denote the categories of presheaves and sheaves on $\bG$. The categories $\Sh\bG$ and $\Pr\bG$ are abelian categories with enough injectives

\subsubsection{}

If $f:G\rightarrow H$ is a morphism of smooth stacks
then $f_*:\Pr \bG\rightarrow \Pr \bH$ is left exact since it is a right adjoint.
We therefore have the right derived functor
$Rf_*:D^+(\Pr_{\Ab}\bG)\rightarrow D^+(\Pr_{\Ab} \bH)$. 

If $g:H\to L$ is a second morphism of smooth stacks, then we have a natural transformation
$$R(g\circ f)_*\to Rg_*\circ Rf_*\ .$$ In fact, let
$F\in D^+(\Pr_{\Ab}\bG)$ be  a lower bounded complex of injective presheaves. Then we choose an injective resolution $f_*F\to J$. Note that $g_*(J)$ represents $Rg_*\circ Rf_*(F)$.
Then  using (\ref{trtz}) the required morphism is defined as the composition
$$R(g\circ f)_*(F)\cong (g\circ f)_*(F)\to g_*\circ f_*(F)\to g_*(J)\cong Rg_*\circ Rf_*(F)\ .$$ 
\begin{lem}\label{rcompo}
If  $f$ is smooth or $g$ is representable, then
$$R(g\circ f)_*\cong Rg_*\circ Rf_*\ .$$
\end{lem}
\proof
If $f$ is smooth, then $f^*$ is exact. 
In this case $f_*$ preserves injectives and we can take $J:=f_*(F)$. 
We can now apply Lemma \ref{comis} in order to see that the natural transformation  $(g\circ f)_*(F)\to g_*\circ f_*(F)$ is an isomorphism.

If $g$ is representable, then $g_*$ is exact by Lemma \ref{zwei11}.
In this case we have again by Lemma \ref{comis} that
$R(g\circ f)_*(F)\cong (g\circ f)_*(F)\cong g_*\circ f_*(F)\cong Rg_*\circ Rf_*(F)$.
\hB

\subsubsection{}

Let $i:\Sh\bG\rightarrow \Pr\bG$ denote the inclusion. It has a left adjoint
$i^{\sharp}:\Pr\bG\rightarrow \Sh\bG$, the sheafification functor (see \cite[3.1.1, 3.2.1]{MR1317816}).
Since the functor $i$ is a right adjoint, it is left exact. We can form
its right derived $Ri:D^+(\Sh_{\Ab}\bG)\rightarrow D^+(\Pr_{\Ab}\bG)$.

Let $f:G\to H$ be a morphism of smooth stacks.
\begin{lem}\label{ddddw}
The functor $i$ preserves injectives and 
we have an isomorphism
$R(f_*\circ i)\cong Rf_*\circ Ri$.
\end{lem}
\proof
Since $i^\sharp$ is exact (see \cite[Thm. 3.2.1 (ii)]{MR1317816}) the functor $i$ preserves injectives.
This implies the assertion.
\hB

\subsubsection{}

Let $\tau:=(U_i\to U)_{i\in I}\in \cov_\bG(U\to G)$ be a covering family of $(U\to G)\in \bG$. For a presheaf $F\in \Pr_\Ab\bG$ 
we form the \v{C}ech complex $\check{C}^*(\tau,F)$. Its $p$th group  is  
$$\check{C}^p(\tau,F):=\prod_{(i_0,\dots,i_p)\in I^{p+1}} F(U_{i_0}\times_U\dots\times_UU_{i_p})\ ,$$ and the differential is given by the usual formula.

\begin{ddd}[see  3.5.1, \cite{MR1317816}]\label{flabbydef}
A sheaf $F\in \Sh_\Ab \bG$ is called flabby if for all $(U\to G)\in \bG$ and all $\tau\in \cov_{\bG}(U\to G)$ we have
$H^k(\check{C}(\tau,F))\cong 0$ for all $k\ge 1$.
\end{ddd}
 
\subsubsection{}

Let $f:\bG\to \bH$ be a smooth map between smooth stacks.
\begin{lem}\label{flabbypres}
The functor $f^*:\Pr_\Ab\bG\to \Pr_\Ab\bH$ preserves flabby sheaves.
\end{lem}
\proof
We have the functor $f_\sharp:\bG\to \bH$ given by $f_\sharp(V\to G):=(V\to G\to H)$ (see \ref{gg221}).
By Lemma \ref{poi451} we know that $f^*$ preserves sheaves.
 
Let $(U\to G)\in \bG$ and $\tau:=(U_i\to U)\in \cov_\bG(U)$.
Observe that $f_\sharp\tau:=(f_\sharp(U_i)\to f_\sharp(U))$ is a covering family of $f_\sharp U$ in $\bH$.

Let $F\in \Sh_\Ab\bH$. By Lemma \ref{ttgg221} we have  $f^*F(U)\cong F(f_\sharp U)$.
We therefore have an isomorphism of complexes
$$\check{C}^\cdot(\tau,f^*F)\cong \check{C}^\cdot(f_\sharp\tau,F)\ .$$ 
If $F$ is in addition flabby, then the
cohomology groups of the right-hand side in degree $\ge 1$ vanish.

\hB

\subsubsection{}

Let $f:G\to H$ be a representable map between smooth stacks.  
\begin{lem}\label{klom}
If $F\in \Sh_\Ab\bG$ is flabby, then $F$ is $(f_*\circ i)$-acyclic.
\end{lem}
\proof
Let $F\in \Sh_\Ab\bG$ be flabby.
We must show that
$R^k(f_*\circ i)(F)\cong 0$ for $k\ge 1$.
By Lemma \ref{zwei11} the functor
$f_*$ is exact so that $R^k(f_*\circ i)(F)\cong f_*\circ R^ki(F)$. Since $F$ is injective it is flabby. Since flabby sheaves are $i$-acyclic
by \cite[Corollary  3.5.3]{MR1317816} we get $R^ki(F)\cong 0$.
\hB

\subsubsection{}

Let $G$ be a smooth stack, $F\in D^+(\Sh_{\Ab}\bG)$,
  and $p:G\rightarrow *$ the canonical morphism. Then we have the object
$Rp_* \circ Ri\in D^+(\Pr_{\Ab} \Site(*))$. Let $\ev:\Pr_{\Ab}\Site(*)\rightarrow \Ab$ be the
evaluation at the object $(*\rightarrow *)\in \Site(*)$. This functor is exact.
\begin{ddd}\label{cogo}
We define the cohomology of $F\in D^+(\Sh_{\Ab}\bG)$ as
$$h(G;F):=\ev\circ Rp_* \circ Ri(F)\in D^+(\Ab)$$
Furthermore we set
$H^*(G;F):=H^*h(G;F)$.
\end{ddd}

In particular, for an abelian group $Z$ we have the constant
presheaf $Z_G$ with value $Z$.
\begin{ddd} We define the cohomology of the smooth stack $G$ with coefficients in $Z$ by 
$$H^*(G;Z):=H^*(G;i^\sharp Z_G)\ .$$
\end{ddd}

\subsubsection{} \label{site-comp}
In \cite[p. 19/20]{math.DG/0605694} another site is used
for sheaves on a smooth stack and
their (hyper)cohomology. In the language of \cite{math.DG/0605694}
a stack is represented as a fibered category over $\Mf^\infty$,
and the open covering topology is used on the underlying
category to define sheaves and cohomology.
This site is equivalent to the site $\Site^a(G)$
of {\em arbitrary}
maps from smooth manifolds to the stack $G$ equipped
with the open covering topology which contains more
objects than $\Site(G)$. In \cite{heinloth} also
the site $\Site^a(G)$ is used.
We have the embedding
$\varphi_G \colon \Site(G) \to \Site^a(G)$
which gives rise to an exact restriction functor
$\varphi_G^* \colon \Sh_{\Ab} \Site^a(G) \to
\Sh_{\Ab} \Site(G)$.
The cohomology $h(G;F)$ can also be defined as
the right derivation of the global sections
functor $\Gamma \colon \Sh_{\Ab} \Site(G) \to \Ab$.
In \cite{math.DG/0605694} the cohomology is defined as
the right derivation of the analogous global
sections functor $\Gamma^a \colon \Sh_{\Ab} \Site^a(G) \to \Ab$.
By universality and the fact that
global sections commute with the
restriction $\varphi_G^*$ there is an induced transformation
$R \Gamma^a \to R \Gamma \circ R \varphi_G^*$.
One shows that this is an isomorphism by using
that $\varphi_G^*$ preserves flabby sheaves,
and the simplicial model description of
the cohomology of Section \ref{toola} which works
for both sites, and is used in \cite{MR2172499} as well as in the present
paper.

\subsection{Simplicial models}\label{toola}

\subsubsection{}

For a morphism $f:G\to X$ between smooth stacks we defined a functor $f_*:\Pr\bG\to \Pr \bX$ (see Definition \ref{dd3}).  We are in particular interested in its derived version 
$Rf_*\circ Ri:D^+(\Sh_\Ab \bG)\to D^+(\Pr_\Ab \bX)$. The definitions of $f_*$ in terms of a limit, and of $Rf_*$  using injective resolutions are very useful for
the study of the functorial properties of $f_*$.  For explicit calculations we would like to work with more concrete objects. In the present  subsection we associate to a flabby sheaf $F\in \Sh_\Ab\bG$ an explicit complex of presheaves $C^\cdot_A(F)\in C^+(\Pr_\Ab\bX)$ which represents $Rf_*\circ i(F)\in D^+(\Pr_\Ab\bX)$ (see Lemma \ref{schluss}). 
It looks like a presheaf of \v{C}ech complexes and depends on the choice of a surjective smooth and representable map  $A\to G$ such that $A\to G\to X$ is also representable (e.g. an atlas of $G$). 

In the present paper we consider three applications of this construction.
The first is the derived version of Lemma \ref{pulb} which says that pull-back and push-forward in certain cartesian diagrams commute (see Lemma \ref{derpulb}).
In the second application we use the complex $C_A^\cdot$ in order to get a de Rham model of the derived push-forward of the constant sheaf with value $\R$
on $\bG$ (see equation (\ref{derham-model}). Finally we use this construction in Lemma \ref{bas13} in  order to calculate the cohomology of the gerbe $[*/S^1]$ explicitely.

\subsubsection{}
Let $G$ be a smooth stack and $(A\to G), (B\to G)\in \bG$.
\begin{lem}\label{fibcat}
The fibre product in stacks
$$\xymatrix{&A\times_GB\ar[dd]\ar[dr]\ar[dl]&\\A\ar[dr]&\ar@{=>}[l]^u&B\ar[dl]\ar@{=>}[l]^{\id}\\&G&}$$
is the categorical product $(A\to G)\times_{\bG}(B\to G)$.
\end{lem}
\proof
The fibre product $(H\to G),(L\to G)\mapsto H\times_GL$ of stacks $H,L\in \cC$ over $G$
is the two-categorical fibre product in the two-category $\cC/G$ of stacks over $G$.
Let $\cC_0\subset\cC$ be the full subcategory of stacks which are equivalent to smooth manifolds, i.e. the essential image of the Yoneda embedding $\Mf^\infty\to \cC$. We define the one-category
$\overline{\cC_0/G}$ by identifying two-isomorphic morphisms and observe that the canonical functor
$\cC_0/G\to \overline{\cC_0/G}$ is an equivalence. Under this equivalence the
restriction of the fibre product to $\cC_0$ becomes the one-categorical product.
This implies the result since the natural functor
$$\bG\to \cC_0/G\to \overline{\cC_0/G}$$ 
is an equivalence of categories. \hB

% Alternatively we can
% verify the universal property.
% To this end we consider the diagram
% $$\xymatrix{&A\times_GB\ar[ddl]\ar[ddr]&\\&T\ar@{.>}[u]^s\ar[dd]\ar[dr]\ar[dl]&\\A\ar[dr]&\ar@{=>}[l]&B\ar[dl]\ar@{=>}[l]\\&G&}\ .$$
% We get the canonical map $s$ and a canonical two-morphism which determines the canonical morphism $(T\to G)\to (A\times_G B\to G)$ in $\bG$ such that all required diagrams two-commute in the correct way. We leave this diagram chase to the interested reader. \hB 
%  

\subsubsection{}\label{gensi}

Let $f:G\rightarrow X$ be a map of smooth stacks.
Let further $A$ be a smooth stack and $A\rightarrow G$ be a representable, surjective and smooth map such that the composition $A\to G\to X$ is also representable. An atlas of $G$ would have these properties, but in applications we will need this more general situation where $A$ is not necessarily  equivalent to a manifold.
Let $(U\to X)\in \bX$ and  form the following diagram of cartesian squares:
\begin{equation}\label{dopw123}\xymatrix{
A_U\ar[d]\ar[r]&A\ar[d]\\
G_U\ar@{=>}[ur]\ar[r]^{j^U}\ar[d]&G\ar[d]\\ 
U\ar@{=>}[ur]\ar[r]&X}\ .
\end{equation}
Since smoothness is preserved by pull-back
the horizontal maps are smooth.
Since surjectivity is also preserved by pull-back the two upper vertical maps are surjective and smooth. Since $A\to X$ is representable, the stack $A_U\cong U\times_XA$ is equivalent to a manifold.

\subsubsection{}
\newcommand{\bGU}{\mathbf{G_U}}

Note that $(A_U\to G_U)\in \bGU$.
In view of Lemma \ref{fibcat} we can take powers of $A_U$ in  $\bGU$.
Using these powers we form a simplical object
$A_U^\cdot\in \bGU$. Its $n$-th object is given by  
$$\underbrace{A_U\times_{G_U}\dots\times_{G_U}A_U}_{n+1\: factors}\to G_U\ .$$
We let $j^U_!A_U^\cdot\in \bG$ denote the simplicial object in $\bG$ with $n$th object
$(A_U^n\to G_U\stackrel{j^U}{\to} G)$.
If $V\to U$ is a morphism in $\bX$, then we obtain an induced morphism of simplicial objects
$j^V_!A^\cdot_V\to j^U_!A^\cdot_U$ in $\bG$.

% We obtain a groupoid in $\bG$
% $$\xymatrix{&A\times_GA\ar[dd]\ar[dr]\ar[dl]&\\A\ar[dr]_s&\ar@{=>}[l]^u&A\ar[dl]^s\ar@{=>}[l]^{\id}\\&G&}\ .$$  which extends to a simplicial object $A^\cdot\in \bG$ such that
% $A^n=(\underbrace{A\times_G\dots\times_GA}_{n+1\times}\rightarrow G)$ and the map to $G$ is induced from the last factor. 

\subsubsection{}\label{check}\label{localo}

If $F\in \Pr\bG$, then we consider the cosimplicial object
$U\mapsto F(j^U_!A_U^\cdot)$ in $\Pr\bX$.  For a morphism $V\to U$ in $\bG$
the structure map $F(j^U_! A_U^\cdot)\to F(j^V_!A_V^\cdot)$ is induced by the morphism of simplicial objects $j^V_!A^\cdot_V\to j^U_!A^\cdot_U$ in $\bG$. 

\begin{ddd}
For a presheaf of abelian groups $F\in {\Pr}_\Ab \bG$ let
$$C_A^\cdot(F)\in C^+({\Pr}_\Ab \bX)$$ denote the chain complex of presheaves associated to the cosimplicial presheaf of abelian groups $U\mapsto F(j^U_!A_U^\cdot)$.
Its differential will be denoted by $\delta$. 
\end{ddd}

\subsubsection{}

Let $F\in \Pr_\Ab\bG$.
\begin{lem}\label{zz}
We have a natural transformation $\psi:f_*F\rightarrow H^0C_A^\cdot(F)$ which is an isomorphism if $F$ is a sheaf.
\end{lem}
\proof
Let $(U\to X)\in \bX$.
We recall definition of the push-forward \ref{pushdef}: 
$$f_*F(U)=\lim_{(V\to G)\in\bG/U}  F(V)\ .$$
Observe that
$$
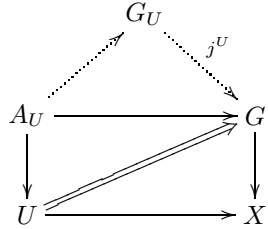
\xymatrix{&G_U\ar@{.>}[dr]^{j^U}&\\A_U\ar@{.>}[ur]\ar[d]\ar[rr]&&G\ar[d]\\
U\ar[rr]\ar@{=>}[urr]&&X}$$
belongs to $\bG/U$  so that we have an evaluation
$$\xymatrix{f_*F(U)\ar[rr]^{evaluation}\ar[dr]^\psi&&C^0_A(F)(U)\\&H^0C_A^\cdot(F)(U)\ar[ur]&}$$
with a canonical factorization $\psi$ by the definition of $H^0C_A^\cdot(F)(U)$ as a kernel.

Assume now that $F$ is a sheaf. Then we must show that $\psi$ is an isomorphism.
Let 
\begin{equation}\label{uoverg1}\xymatrix{V\ar[d]\ar[r]&G\ar[d]\\U\ar@{=>}[ur]\ar[r]&X}\end{equation}
 be in  $\bG/U$.
Then we have a canonical factorization
$$
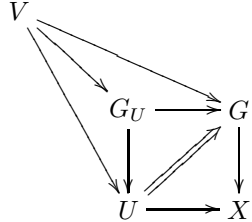
\xymatrix{V\ar[dr]\ar[ddr]\ar[drr]&&\\&G_U\ar[d]\ar[r]&G\ar[d]\\&U\ar@{=>}[ur]\ar[r]&X}\ .$$
Using the induced map $V\to G_U$ we form the diagram
$$\xymatrix{V\times_{G_U}A_U\ar[r]\ar[d]&A_U\ar[d]\ar[dr]&\\
V\ar@{=>}[ur]\ar[r]&G_U\ar[r]^{j^U}&G}$$
We consider the composition $(V\times_{G_U}A_U\to A_U\to G)$ as an object in $\bG$.  Since $A_U\to G_U$ is smooth and surjective the map
$(V\times_{G_U}A_U\rightarrow V)$ is a covering of $V$ in $\bG$ (it is here where we use the submersion-pre-topology). For a sheaf $F$ we have
$$F(V)\cong \lim\left(\xymatrix{F(A_U\times_{G_U}V)\ar@{=>}[r]&F((A_U\times_{G_U} V)\times_V(A_U\times_{G_U} V))}\right)\ .$$
We further have 
$$(A_U\times_{G_U} V)\times_V(A_U\times_{G_U} V)\cong A_U\times_{G_U}A_U\times_{G_U}V$$ and a diagram
$$\xymatrix{F(A_U\times_{G_U}V)\ar@{=>}[r]&F(A_U\times_{G_U}A_U\times_{G_U}V)\\F(j_!^UA^0_U)\ar@{=>}[r]\ar[u]&F(j_!^UA^1_U) \ar[u]}$$
(recall that $A_U^0=A_U$ and $A_U^1=A_U\times_{G_U}A_U$)
induced by the projection along $V$.
Since $H^0C_A^\cdot(F)(U)$ is the limit of the lower horizontal part
the left vertical map  induces a map
$H^0C_A^\cdot(F)(U)\rightarrow F(V)$. 
Since this construction is natural in the object (\ref{uoverg1}) of $\bG/U$ we obtain finally a map
$H^0C_A^\cdot(F)(U)\rightarrow f_*F(U)$ which is the inverse to $\psi$.
\hB 

\subsubsection{}

\begin{lem}\label{injvan}
If $F\in \Pr_\Ab\bG$ is injective, then
$H^iC_A^\cdot(F)=0$ for $i\ge 1$.
\end{lem}
\proof
We follow the ideas of the last part of the  proof of  \cite[Thm. 2.2.3]{MR1317816}.
Let $(U\to X)\in \bX$ and
$\cA_U^\cdot$ denote the simplicial presheaf of sets represented by
$j^U_!A_U^\cdot$. Furthermore, let $\Z_{\cA_U^\cdot}$ be the (non-positively  graded) complex of free abelian presheaves generated by $\cA_U^\cdot$.
Then for any presheaf $F\in \Pr_\Ab \bG$ we have
$$C_A^\cdot(F)(U)\cong \Hom_{\Pr_\Ab\bG}(\Z_{\cA_U^\cdot},F)\ .$$

Since $F$ is injective $\Hom_{\Pr_\Ab\bG}(\dots ,F)$ is an exact functor.  Hence it suffices to show that $H^i(\Z_{\cA_U^\cdot})=0$ for $i\le -1$. 
For $(V\to G)\in \bG$ the complex
$\Z_{\cA_U^\cdot}(V)$ is the complex associated to the linearization of the simplicial set
$\Hom_{\bG}(V,j^U_!A_U^\cdot)$.
We now rewrite
\begin{equation}\label{onecatego}A_U\times_{G_U}\dots \times_{G_U}A_U\cong (A\times_G\dots\times_GA)\times_XU\cong (A\times_G\dots\times_GA)\times_GG_U\ .
\end{equation}
We consider $V$, $(A\times_G\dots\times_GA)$ and $G_U$ with their canonical maps to $G$  as objects of the two-category $\cC/G$ of stacks over $G$.  The first object is a manifold and therefore
does not have non-trivial two-automorphisms. Since the maps $(A\times_G\dots\times_GA)\to G$ and  $G_U\to G$ are representable these objects of  $\cC/G$ also do not have non-trivial two-automorphisms. By the same reasoning as in the proof of Lemma \ref{fibcat} we can interpret the fibre product (\ref{onecatego})
as a one-categorical product.  We get
\begin{eqnarray*}
\Hom_{\bG}(V,j^U_!A_U^\cdot)&=&\Hom_{\cC/G}(V,(A\times_G\dots\times_GA)\times_GG_U)\\
&\cong&(\Hom_{\cC/G}(V,A)\times \dots\times \Hom_{\cC/G}(V,A))\times \Hom_{\cC/G}(V,G_U) \\
&\cong&\Hom_{\cC/G}(V,A)^\cdot\times \Hom_{\cC/G}(V,G_U) 
\end{eqnarray*}
For any set $S$, if we take the simplicial set $S^\cdot$ of the powers of $S$, the complex associated to the linearization $\Z_{S^\cdot}$ is exact in degrees $\le -1$. 
Therefore  the complex
$\Z_{\Hom_{\cC/G}(V,A)^\cdot}$ is exact in degree $\le -1$.
Since the tensor product with the free abelian group  $\Z_{ \Hom_{\cC/G}(V,G_U)}$ is an exact functor the complex
$$\Z_{\Hom_{\bG}(V,j^U_!A_U^\cdot)}\cong \Z_{\Hom_{\cC/G}(V,A)^\cdot}\otimes
\Z_{ \Hom_{\cC/G}(V,G_U)} $$
is exact in degree $\le -1$, too.
\hB

\subsubsection{}\label{dpluss}

Since exactness of complexes of presheaves is defined objectwise the functors $C^p_A:\Pr_\Ab \bG\to \Pr_\Ab \bX$ are  exact for all $p\ge 0$. Composing with the total complex construction we  extend the functor $C^\cdot_A$ to a  functor between the categories of lower bounded complexes 
$C_A:C^+(\Pr_{\Ab}\bG)\rightarrow C^+(\Pr_{\Ab}\bX)$ 
(in order to distinguish this from the double complex we drop the ${}^\cdot$
at the symbol $C_A$).
Since this functor is level-wise exact it descends to a functor 
$C_A:D^+(\Pr_{\Ab}\bG)\rightarrow D^+(\Pr_{\Ab}\bX)$
between the lower-bounded derived categories.

\subsubsection{}\label{alg}

Assume that $F$ is a presheaf of associative algebras on $\bG$.
Then $C_A^\cdot(F)$ is a presheaf of $DG$-algebras in the  following natural way.
Pick $$\alpha\in F(j_!^UA_U^{p})\cong C^p_A(F)(U)\ ,\quad \beta\in F(j_!^UA_U^{q})\cong C_A^q(F)(U)\ .$$ 
We have natural maps $u:j_!^UA_U^{p+q}\rightarrow j_!^UA_U^{p}$ and $v:j_!^UA_U^{p+q}\rightarrow j_!^UA_U^{q}$ in $\bG$ projecting onto the first $p+1$ or last $q+1$ factors, respectively. 
Then we define
$\alpha\cdot\beta\in F(j_!^UA_U^{p+q})\cong C_A^{p+q}(F)(U)$ by
$u^*\alpha\cdot v^*\beta$. 
One easily checks that $\delta(\alpha\cdot\beta)=\delta\alpha+(-1)^{p} \alpha\cdot \delta\beta$.

% \subsubsection{}\label{alg1}
% 
% More general, if $F$ is a presheaf of $DG$-algebras, then  $C_A^\cdot(F)$ is again a presheaf of $DG$-algebras. If $F$ is in addition (graded!) commutative, then $H^0C_A^\cdot(F):=\ker(\delta:C_A^0(F)\to C_A^1(F))$ is a presheaf of central $DG$-subalgebras of $C_A(F)$. 

% \subsubsection{}\label{localo}
% 
% We now localize these constructions.
% If $(U\rightarrow X)\in \bX$, then
% we form
% $$\xymatrix{&U\times_XG\ar[d]\ar[r]^v&G\ar[d]^f\\
% A_U:=U\times_XA\ar[ur]^{s_U}\ar[r]&U\ar@{=>}[ur]\ar[r]^u&X}\ .$$
% The map $s_U$ induced by $s$ and smooth. Since $v$ is smooth 
% as a pull-back of the smooth map $u$ the composition $v\circ s_U$ is smooth, too.
% 
% The  prescription  $U\mapsto A_U^\cdot$ is a functor from
% $\bX$ to simplicial objects in $\bG$.
% Let $C^+(\Pr_\Ab \bX)$ denote the category of lower bounded complexes of presheaves of abelian groups on $\bX$.
% We define a functor
% $C_A:\Pr_{\Ab}\bG\rightarrow C^+(\Pr_{\Ab}\bX)$   which maps
% $L\in \Pr_{\Ab}\bG$ to $C_A(L):= (U\mapsto L(A^\cdot_U))$. 

\subsubsection{}\label{shdg}

If $F^\cdot$ is a presheaf of commutative $DG$-algebras,
then $C_A(F^\cdot)$ is a presheaf of associative $DG$-algebras central over the presheaf of commutative $DG$-algebras $(\ker(\delta):C^0_A(F^\cdot)\to C^1_A(F^\cdot))$.

\subsubsection{}

By Lemma \ref{zz} we have a map
$$ \psi:f_*F\rightarrow H^0 C^\cdot_A(F)$$
which is an isomorphism if $F$ is a sheaf.

\begin{lem}
For all  $F\in \Pr_{\Ab}\bG$ we have a natural isomorphism
$RH^0C_A^\cdot(F)\cong C^\cdot_A(F)$ in $D^+(\Pr_{\Ab} X)$.
\end{lem}
\proof
Let $F\rightarrow I^\cdot$ be an injective resolution.
Then we have
$RH^0C^\cdot_A(F)\cong H^0C^\cdot_A(I^\cdot)$.
By  Lemma \ref{injvan} the inclusion
$H^0C_A^\cdot(I^\cdot)\rightarrow C_A(I^\cdot)$ is a quasi-isomorphism.
Since $C_A^\cdot$ is exact the quasi-isomorphism $F\rightarrow I^\cdot$ induces a quasi-isomorphism
$C_A^\cdot(F)\cong C_A(I^\cdot)$. 
\hB

\subsubsection{}

Recall that $i:\Pr_{\Ab}\bG\rightarrow \Sh_{\Ab}\bG$ is left exact and admits a right derived functor
$Ri:D^+(\Sh_{\Ab}\bG)\rightarrow D^+(\Pr_{\Ab}\bG)$, and that
$C_A$ descends to a functor between the lower bounded derived categories (see \ref{dpluss}).
\begin{lem}\label{schluss}
We have a natural isomorphism of functors
$$C_A\circ Ri\cong Rf_*\circ Ri:D^+(\Sh_{\Ab}\bG)\rightarrow D^+({\Pr}_{\Ab}\bX)\ .$$
\end{lem}
\proof
By Lemma \ref{zz} we have an isomorphism of functors
$f_*\circ i\cong H^0C_A\circ i$. 
Hence we have an isomorphism
$$Rf_*\circ Ri\stackrel{!}{\cong} R(f_*\circ i)\cong R(H^0C^\cdot_A\circ i)\stackrel{!}{\cong} RH^0C^\cdot_A\circ Ri\cong   C_A\circ Ri\ ,$$
where at the marked isomorphisms we use that $i$ preserves injectives (compare Lemma \ref{ddddw}).
\hB
 
\subsubsection{}

Assume that we have a diagram in smooth stacks
\begin{equation}\label{squ5671}\xymatrix{G\ar[r]^{u}\ar[d]^f&H\ar[d]^{g}\\
X\ar@{=>}[ur]\ar[r]^v&Y}\ ,\end{equation}
 where $u$  and $v$ are smooth.
Note that $u^*$ and $v^*$ are exact (Lemma \ref{poi451}).

\begin{lem}\label{derpulb}
\begin{enumerate}
\item
We have a natural transformation of functors $D^+(\Pr_\Ab\bH)\rightarrow D^+(\Pr_\Ab\bX)$
$$v^*\circ Rg_*\rightarrow Rf_*\circ u^*\ .$$
\item
The induced transformation $D^+(\Sh_\Ab\bH)\rightarrow D^+(\Pr_\Ab\bX)$
$$v^*\circ Rg_*\circ Ri\rightarrow Rf_*\circ u^*\circ Ri$$
is an isomorphism if
(\ref{squ5671}) is cartesian.
\end{enumerate} 
\end{lem}
\proof
The transformation $(1)$ is induced by
$$v^*\circ Rg_*\cong R(v^*\circ g_*)\stackrel{\ref{pulb}}{\cong} R(f_*\circ u^*)\to Rf_*\circ u^*\ .$$

In order to show the second part $(2)$  we must show that
$$R(f_*\circ u^*)\circ Ri\to Rf_*\circ u^*\circ Ri$$ is an isomorphism.
We calculate $Ri$ using injective resolutions. Note that $i$ preserves injectives. Hence  in order to show that this map is an isomorphism it suffices to show that $u^*$ maps injective sheaves to $f_*$-acyclic presheaves.

Note that $u^*$ preserves sheaves (Lemma \ref{poi451}). We let $u^*_s:\Sh_\Ab \bH\to \Sh_\Ab \bG$ denote the restriction of $u^*$ to sheaves.
Let $F\in \Sh_\Ab\bH$ be injective.  Since injective sheaves are flabby, flabby sheaves are $i$-acyclic,  and $u^*$ preserves flabby sheaves (see Lemma \ref{flabbypres}) we have
$$Rf_*\circ i\circ u_s^*(F)\cong Rf_*\circ Ri\circ u^*_s(F)\stackrel{Lemma\:\ref{schluss}}{\cong}  C_A\circ Ri\circ u^*_s(F)\cong C_A\circ u^*\circ i(F)\ .$$

We now show that the higher cohomology presheaves of $C^\cdot_A\circ u^*\circ i(F)$ vanish. Let $(U\to X)\in \bX$ and choose an atlas
$B\to H$. Then we get the following extension of the diagram (\ref{dopw123})
\begin{equation}\label{dopw123ex}\xymatrix{
A_U\ar[d]\ar[r]&A\ar[d]\ar[r]&B\ar[d]\\
G_U\ar@{=>}[ur]\ar[r]^{j^U}\ar[d]&G\ar@{=>}[ur]\ar[d]\ar[r]^u&H\ar[d]^g\\ 
U\ar@{=>}[ur]\ar[r]&X\ar@{=>}[ur]\ar[r]^v&Y}
\end{equation}
such that all squares are cartesian. The three upper vertical maps are smooth and surjective.  The composition $A\to G\to X$ is representable. All horizontal maps are smooth.
We have the simplicial object $(A_U^\cdot\to G_U)\in \bGU$ and let
$u_!j^U_! A_U^\cdot\in \bH$ be the induced simplicial object
$A_U^\cdot\to G_U\stackrel{j^U}{\to} G\stackrel{u}{\to} H$ in $\bH$.
Then we have by the definition \ref{check} of $C_A$  and the formula (\ref{ttgg221}) for $u^*$ that
$$C_A^\cdot\circ u^*\circ i(F)(U)=F(u_!j^U_! A_U^\cdot)\ .$$ 
We now observe the isomorphisms
\begin{eqnarray*}
A_U\times_{G_U}\dots\times_{G_U}A_U&\cong&(A\times_G\dots\times_GA)\times_XU\\
&\cong&(B\times_H\dots\times_HB)\times_YU\\
&\cong&B_{v_!U}\times_{H_{v_!U}}\dots\times_{H_{v_!U}}B_{v_!U}\ ,
\end{eqnarray*}
where the notation is explained by the cartesian diagram
$$\xymatrix{H_{v_!U}\ar[d]\ar[r]^k&H\ar[d]\\U\ar[r]&Y}\ ,$$
and where
$v_!U:=(U\to X\stackrel{v}{\to} Y)\in \bY$. We can thus identify
the simplicial object $u_!j_!A_U^\cdot$ with the similar simplicial object 
$k_!B_{v_! U}$ in $\bH$. In other words, we have an isomorphism of complexes
$$C_A^\cdot\circ u^*\circ i(F)(U)\cong C_B^\cdot\circ i(F)(v_!U)\ .$$
Since $i(F)$ is an injective presheaf the right-hand side is exact by Lemma \ref{injvan}. \hB

\subsection{Comparison of big and small sites}\label{cckkll11}

\subsubsection{}

Let $X$ be a smooth stack and $(U\to X)\in \bX$. A presheaf on $X$ naturally induces a presheaf on the small site $(U)$ of the manifold $U$ consisting of the open subsets. This restriction functor will be used subsequently in order to reduce assertions in the sheaf theory over $X$ to
assertions in the ordinary sheaf theory on $U$. 
The goal of the present subsection is to study exactness properties of this restriction
and its relation with the sheafification functors.

\subsubsection{}\label{cdsa1}

If $U$ is a smooth manifold, then we let $(U)$ denote the small site of
$U$ where covering families are coverings by families of open submanifolds. 
A presheaf with respect to the big site on $U$ is in particular a presheaf with respect to
$(U)$.

\subsubsection{}

Let $G$ be a smooth stack and $(U\to G)\in \bG$.
Then we have a functor $\nu_U:\Pr\bG\to \Pr(U)$ which associates to the presheaf $F\in \Pr\bG$ the presheaf $\nu_U(F)\in \Pr(U)$ obtained by restriction of structure.  Since limits and colimits in presheaves are defined objectwise the functor $\nu_U$ is exact.

\subsubsection{}

\begin{lem}\label{sheafres1}
The functor $\nu_U$ preserves sheaves and induces a functor $\nu_U^s:\Sh\bG\to \Sh(U)$. 
\end{lem}
\proof
An object $V\in (U)$ gives rise to an object $(V\to U\to G)\in \bG$. 
Observe that covering families of objects of $V\in (U)$ are also covering families of $(V\to G)\in \bG$.
For open subsets $V_1,V_2\subset V$ the fibre products $V_1\times_V V_2$ in
$(U)$ and in $\bG$ coincide by the discussion in \ref{t5tt5t}.
Therefore the descent conditions on $\nu_U(F)$ to be a sheaf on $(U)$ are part of the descent conditions for $F$ to be a sheaf on $\bG$.  
Hence the functor $\nu_U$ restricts to
$\nu_U^s:\Sh\bG\to \Sh(U)$.

\hB

\subsubsection{}
 Since limits of sheaves are defined objectwise the functor $\nu_U^s$ commutes with limits. The goal of the following discussion is to show that it also commutes with colimits.
\begin{prop}\label{exgt5}
The functor $\nu_U^s:\Sh\bG\to \Sh(U)$ is exact.
\end{prop}
\proof 
If $F$ is a diagram of sheaves, then we have
$$\colim^s(F)\cong i^\sharp\circ \colim\circ i(F)\ ,$$
where $\colim^s$ is the colimit of sheaves.  Note that $\nu_U\circ i\cong i\circ \nu^s_U$ and $\nu_U\circ \colim\cong \colim\circ \nu_U$. In order to show that $\nu^s_U$ commutes with $\colim^s$ it remains to show the following Lemma. 

\subsubsection{}
\begin{lem}\label{gloloc}
We have
$$i^\sharp\circ \nu_U\cong \nu_U^s\circ i^\sharp:\Pr\bG\to \Sh(U)\ .$$
\end{lem}
\proof
For the moment it is useful to indicate by a subscript (e.g. $i_\bG$ or $i_{(U)}$) the site for which the functors are considered.
Following the discussion in \cite[Section 3.1]{MR1317816} we introduce
an explicit construction of the sheafification functor.
Consider the site $\bG$.
We define the functor $P_\bG:\Pr\bG\to \Pr\bG$ as follows.
Let $(V\to G)\in \bG$. Then we have the category of covering families $\cov_\bG(V)$ whose morphisms are refinements. For $\tau:=(V_i\to V)\in \cov_\bG(V)$ we define
$H^0(F)(\tau)$ by the equalizer diagram
$$H^0(F)(\tau)\to\prod_{i} F(V_i)\Longrightarrow \prod_{i,j}F(V_i\times_V V_j)\ .$$
We get a diagram
$\tau\to H^0(F)(\tau)$ in $\Sets^{\cov_\bG(V)}$ and define
$$P_\bG(F)(V):=\colim_{\tau\in \cov_\bG(V)} H^0(F)(\tau)\ .$$
Then we have $$i_\bG\circ i_\bG^\sharp:=P_\bG\circ P_\bG:\Pr\bG\to \Pr\bG\ .$$

In a similar manner we define a functor $P_{(U)}:\Pr(U)\to \Pr(U)$ and get
$$i_{(U)}\circ i_{(U)}^\sharp:= P_{(U)}\circ P_{(U)}:\Pr(U)\to \Pr(U)\ .$$

In order to show the Lemma it suffices to show that
$$P_{(U)}\circ \nu_U\cong \nu_U\circ P_\bG\ .$$

Let $V\subset U$ be open and consider the induced $(V\to
G)\in \bG$. Then we have a functor
\begin{equation}\label{cifin}a:\cov_{(U)}(V)\to \cov_{\bG}(V)\ .\end{equation} 
If $\tau\in \cov_{(U)}(V)$, then we have an isomorphism
$$H^0(F)(a(\tau))\cong H^0(\nu_U(F))(\tau)\ .$$
We therefore have an induced map of colimits
$$P_{(U)}\circ \nu_U(F)(V)\to \nu_U\circ P_{\bG}(F)(V)\ .$$
This map is in fact an isomorphism since we will show below that
(\ref{cifin}) defines a cofinal subfamily.

Let $\sigma:=(U_i\to V)_{i\in I}\in \cov_{\bG}(V)$. Since the maps $U_i\to V$ are submersions they admit local sections. Hence there exists a covering $\tau:(V_j\to V)_{i\in J}\in \cov_{(U)}$, a map
$r:J\to I$ and a family of sections $s_j:V_j\to U_{r(j)}$ such that
$$\xymatrix{&U_{r(j)}\ar[dr]\\V_j\ar[ur]^{s_j}\ar[rr]&&V}
$$ commutes for all $j\in J$. 
This data defines a morphism $\sigma\to a(\tau)$ in $\cov_\bG(V)$. 
\hB 
This finishes the proof of Proposition \ref{exgt5}. \hB

\subsubsection{}
 Recall the definition of a flabby sheaf \ref{flabbydef}.

\begin{lem}\label{flow12}
The functor $\nu_U^s:\Sh_\Ab\bG\to \Sh_\Ab(U)$ preserves flabby sheaves.
\end{lem}
\proof Let $V\subset U$ be an open subset and  $\tau\in \cov_{(U)}(V)$. 
Let $a:\cov_{(U)}(V)\to \cov_\bG(V)$ be
as in (\ref{cifin}). We have an natural isomorphism $\check{C}(\tau,\nu^s_U(F))\cong \check{C}(a(\tau),F)$. If $F\in \Sh_\Ab\bG$ is flabby, then for $k\ge 1$ we have
$H^k\check{C}(\tau,\nu^s_U(F))\cong H^k(\check{C}(a(\tau),F))
\cong 0$.
\hB

\newcommand{\Mod}{{\tt Mod}}

\subsubsection{}\label{derdecs1}

Since $\nu_U:\Pr\bG\to \Pr (U)$ is exact it descends to a functor 
$\nu_U: D^+(\Pr_\Ab\bG)\to D^+(\Pr_\Ab(U))$ between the lower bounded derived categories. Since $\nu_U^s:\Sh\bG\to \Sh(U)$ is exact, it descends to a functor
$\nu^s_U: D^+(\Sh_\Ab\bX)\to D^+(\Sh_\Ab(U))$.

\subsubsection{}
Using the techniques above we show the following result which will be useful later.
Let $f:G\to H$ be a smooth map between smooth stacks. Note that $f^*:\Pr_\Ab\bH\to \Pr_\Ab\bG$ is exact.
\begin{lem}\label{somecomu}
\begin{enumerate}
\item We have an isomorphism of functors
$$f^*\circ i_\bH\circ i_\bH^\sharp\cong i_\bG\circ i_\bG^\sharp\circ f^*:{\Pr}_\Ab \bH\to {\Pr}_\Ab \bG\ .$$
\item
We have an isomorphism of functors
$$f^*\circ Ri_\bH\circ i_\bH^\sharp\cong Ri_\bG\circ i_\bG^\sharp\circ f^*:D^+({\Pr}_\Ab \bH)\to D^+({\Pr}_\Ab \bG)\ .$$
\end{enumerate}
\end{lem}
\proof
Let $(U\to G)\in \bG$ and $f_!U:=(U\to G\stackrel{f}{\to} H)\in \bH$.
We calculate for $F\in \Pr_\Ab \bH$ that on the one hand
\begin{eqnarray*}
(f^*\circ i_\bH\circ i_\bH^\sharp F)(U)&\stackrel{(\ref{poi451})}{\cong}&(i_\bH\circ i_\bH^\sharp F)(f_{!}U)\\&\cong&(\nu_{f_!U}\circ i_\bH\circ i_\bH^\sharp F)(U)\\
&\stackrel{\ref{sheafres1} ,\ref{gloloc}}{\cong}&(i_{(U)}\circ i_{(U)}^\sharp\circ  \nu_{f_!U}F)(U)\ .
\end{eqnarray*}
On the other hand we have
\begin{eqnarray*}
(i_\bG\circ i_\bG^\sharp\circ f^*F)(U)&\cong &(\nu_U\circ i_\bG\circ i_\bG^\sharp\circ f^*F)(U)\\
&\stackrel{\ref{sheafres1} ,\ref{gloloc}}{\cong}&(i_{(U)}\circ i_{(U)}^\sharp\circ \nu_U \circ  f^*F)(U)\ .
\end{eqnarray*}
Finally we use the fact that $\nu_U\circ f^*F\cong \nu_{f_!U} F$.  Indeed, for
$V\subset U$ we have $$\nu_U\circ f^*F(V)\cong  F(f_!V)\cong  \nu_{f_!U} F(V)\ .$$
The combination of these isomorphisms gives the first assertion.

Since $f^*$ preserves sheaves we can consider the restriction $f_s^*:\Sh_\Ab\bH\to \Sh_\Ab \bG$ of $f^*$. Using the first part of the Lemma  and the isomorphism $i^\sharp\circ i\cong \id$ we get
\begin{equation}\label{teiluo001}f_s^*\circ i^\sharp_\bH\cong i^\sharp_\bG\circ i_\bG\circ f_s^*\circ i_\bH^\sharp\cong i_\bG^\sharp\circ f^*\circ i_\bH\circ i_\bH^\sharp\stackrel{(1)}{\cong}
i_\bG^\sharp\circ i_\bG\circ i_\bG^\sharp\circ f^*\cong i_\bG^\sharp\circ f^*\ .
\end{equation}
Note that $f_s^*$ is an exact functor. In order to see that it is right exact we use that  $f_*$ preserves sheaves we consider its restriction $f_*^s$ to sheaves. For $X\in \Sh\bH$ and $Y\in \Sh \bG$ we get a natural isomorphism
\begin{eqnarray*}
\Hom_{\Sh\bG}(f_s^*X,Y)&\cong&\Hom_{\Sh\bG}(i^\sharp_\bG \circ i_\bG f_s^*X,Y)\\
&\cong&\Hom_{\Pr\bG}( i_\bG f_s^*X, i_\bG Y)\\
&\cong&\Hom_{\Pr\bG}(f^*\circ i_\bH X, i_\bG Y)\\
&\cong&\Hom_{\Pr\bH}(i_\bH X, f_*\circ i_\bG Y)\\
&\cong&\Hom_{\Pr\bH}(i_\bH X, i_\bH \circ f^s_* Y)\\
&\cong&\Hom_{\Sh\bH}(X, f^s_* Y)
\end{eqnarray*}
Therefore $f^*_s$ is a left adjoint and therefore right exact.
We now write $$f_s^*\cong i^\sharp_\bG\circ i_\bG\circ  f_s^*\cong i_\bG^\sharp\circ f^*\circ i_\bH\ .$$ Since $i^\sharp_\bG$ is exact, and $f^*$ and $i_\bH$ are left exact (since they are right adjoints, see Lemma \ref{lef} for a left adjoint of $f^*$ ) we conclude that $f^*_s$ is left exact, too

Using the that $f_s^*$ preserves flabby sheaves (Lemma \ref{flabbypres})
and that it is an exact functor
we get
$$f^*\circ Ri_\bH\cong R(f^*\circ i_\bH)\cong R(i_\bG\circ f_s^*)\cong 
Ri_\bG\circ f_s^*\ .$$
Combining this with (\ref{teiluo001}) we get the desired isomorphism
$$f^*\circ Ri_\bH\circ i^\sharp_\bH\cong Ri_\bG\circ f^*_s\circ i^\sharp_\bH\cong 
Ri_\bG\circ i^\sharp_\bG\circ f^*\ .$$
\hB 

\section{The de Rham complex}\label{se2}

\subsection{The de Rham complex is a flabby resolution}

\subsubsection{}

We want to apply Lemma \ref{schluss} to the sheafification $i^\sharp \R_\bG$
of the constant presheaf with value $\R$ on $\bG$. In particular, we  must
calculate $Ri(i^\sharp \R_\bG)$. This can be done by applying $i$ to a flabby
resolution of $i^\sharp \R_\bG$. In the present subsection we introduce the de
Rham complex  $G$ and show that it is a flabby resolution of $i^\sharp
\R_\bG$. The de Rham complex of smooth stacks has also been investigated in
the papers \cite{MR2172499}, \cite[Sec. 3]{math.DG/0605694}, \cite{MR2183389}.

The de Rham complex of $G$ is built from the de Rham complexes of the manifolds $U$ for all $(U\to G)\in \bG$.
 For each $U$ equipped with the topology of the small site it is well known that the de Rham complex resolves the constant sheaf with value $\R$ and
is flabby. 
Our task here  is to extend these properties to the stack $G$  and the big site.

\subsubsection{}

Let $G$ be a smooth stack and fix an integer $p\ge 0$.
We define the presheaf $\Omega^p(G)$ by
$$\Omega^p(G)(U):=\Omega^p(U)\ .$$
If $\phi:U\to V$ is a morphism in $\bG$, then $\Omega^p(G)(\phi):=\phi^*:\Omega^p(V)\to \Omega^p(U)$. Since $\phi^*$ commutes with the de Rham differential  we get a complex $(\Omega^\cdot(G), d_{dR})$ of presheaves.

\begin{lem}\label{resd97}
The presheaf $\Omega^p(G)$ is a sheaf and flabby.
\end{lem}
\proof
 Let $(U\to G)\in \bG$. Observe that
$\nu_U(\Omega^p(G))$ is the presheaf of smooth sections of the vector bundle $\Lambda^pT^*U$. This is actually a sheaf. 
In order to show that
$\Omega^p(G)$ is a sheaf it suffices to show that the unit
$\Omega^p(G)\to i_{\bG}\circ i_{\bG}^\sharp (\Omega^p(G))$ of the adjoint pair $(i_\bG^\sharp,i_\bG)$ is an isomorphism. 
This follows from the calculation 
\begin{eqnarray*}
i_\bG\circ i_\bG^\sharp (\Omega^p(G))(U)&\cong&\nu_U\circ i_\bG\circ i_\bG^\sharp (\Omega^p(G))(U)\\
&\stackrel{Lemma\: \ref{sheafres1}}{\cong}&i_{(U)}\circ \nu_U^s\circ i_\bG^\sharp (\Omega^p(G))(U)\\
&\stackrel{Lemma\: \ref{gloloc}}{\cong}&i_{(U)}\circ i_{(U)}^\sharp\circ \nu_U (\Omega^p(G))(U)\\
&\stackrel{\nu_U (\Omega^p(G))\:is \:a\: sheaf}{\cong}&\nu_U (\Omega^p(G))(U)\ .
\end{eqnarray*}

A sheaf $F\in \Sh_\Ab(U)$ on a paracompact space $U$ is called soft if for all closed subsets $Z\subset U$ the restriction $\Gamma_U(F)\to \Gamma_Z(F)$ is surjective. For a soft sheaf we have  $R^i\Gamma_U(F)\cong 0$ for all $i\ge 1$ (see \cite[Ex. II 5]{MR1299726}). It now follows from \cite[Corollary 3.5.3]{MR1317816}  that a soft sheaf is flabby.

A sheaf of smooth sections of a smooth vector bundle on a smooth manifold is soft. In particular, $\nu_U^s(\Omega^p(G))$ is soft and therefore flabby.

In order to show that the sheaf $\Omega^p(G)$ is flabby  it suffices by \cite[Corollary 3.5.3]{MR1317816} to show
that $R^k i(\Omega^p(G))\cong 0$ for $k\ge 1$.
We calculate 
\begin{eqnarray*}
R^k i(\Omega^p(G))(U)&\cong& \Gamma_U\circ R^k i(\Omega^p(G))\\
&\stackrel{\Gamma_U \:exact}{\cong}&H^k(\Gamma_U\circ Ri(\Omega^p(G)))\\
&\stackrel{H^k\: object wise}{\cong}&H^k(Ri(\Omega^p(G))(U))\\
&\stackrel{definition \:of\:\nu_U }{\cong}&H^k(\nu_U\circ Ri_\bG(\Omega^p(G))(U))\\
&\stackrel{\nu_U\:exact}{\cong}&R^k(\nu_U\circ i_\bG)(\Omega^p(G))(U)\\
&\stackrel{Lemma \:\ref{sheafres1}}{\cong}&R^k(i_{(U)}\circ \nu_U^s)(\Omega^p(G))(U)\\
&\stackrel{Lemma \:\ref{flow12}}{\cong}&R^ki_{(U)}\circ \nu_U^s(\Omega^p(G))(U)\\
&\stackrel{\nu_U^s(\Omega^p(G))\:is\:flabby}{\cong}&0 \ . 
\end{eqnarray*}
\hB

\subsubsection{}\label{deh}

% By \cite[Prop. 3.4.2]{MR1317816} for $F\in \Sh_{\Ab}(\bG)$ and $(U\rightarrow G)$ in $\bG$ we have
% $R^qi(F)(U)\cong R^q\Gamma_U(F)$, where $\Gamma_U:\Sh_{\Ab}(\bG)\rightarrow \Ab$ is the section functor of $U$. If $F\rightarrow S^\cdot$ is a resolution such that $R^i\Gamma_U(S^p)\cong 0$ for $i\ge 1$ and $p\in \Z$, then we have a natural isomorphism of
% $Ri(F)\cong i(S^\cdot)$ in $D^+(\Pr_{\Ab}\bG)$.

Let $\R_\bG$ denote the constant presheaf on $\bG$  with value $\R$ and $i^\sharp \R_\bG$ its sheafification.
We have a canonical map $\R_\bG\to  H^0(\Omega^\cdot(G))$.
Since $ H^0(\Omega^\cdot(G))$ is a sheaf we get an induced
map
$i^\sharp \R_G\to H^0(\Omega^\cdot(G))$.
\begin{lem}\label{red98}
The map
$$i^\sharp \R_G\to \Omega^\cdot(G)$$ is a quasi-isomorphism.
\end{lem}
\proof
Note that the cohomology sheaves $H^k_s(F^\cdot)\in \Sh_\Ab\bG$ of a complex $F^\cdot$ of sheaves of abelian groups on $\bG$ are defined by $H^k_s(F^\cdot):=i^\sharp_\bG\circ H^k\circ i_\bG(F^\cdot)$, where $H^k$ takes the cohomology of a complex of presheaves objectwise.
We calculate
\begin{eqnarray*}
H_s^k(\Omega^\cdot(G))(U)&\stackrel{definition\:of\: H^k_s}{\cong}&(i_\bG^\sharp\circ H^k\circ i_\bG)(\Omega^\cdot(G))(U)\\
&\stackrel{Lemma \:\ref{sheafres1}}{\cong}&(\nu_U^s\circ i_\bG^\sharp\circ H^k\circ i_\bG)(\Omega^\cdot(G))(U)\\ 
&\stackrel{Lemma  \ref{gloloc}}{\cong}&(i_{(U)}^\sharp\circ \nu_U\circ H^k\circ i_\bG)(\Omega^\cdot(G))(U)\\
&\stackrel{\nu_U\: is\: exact}{\cong}&(i_{(U)}^\sharp \circ H^k \circ \nu_U \circ i_\bG)(\Omega^\cdot(G))(U)\\
&\stackrel{Lemma \:\ref{sheafres1}}{\cong} &(i_{(U)}^\sharp \circ H^k \circ i_{(U)}\circ \nu^s_U)(\Omega^\cdot(G))(U)\\
&\stackrel{definition\:of\: H^k_s}{\cong}&(H_s^k\circ \nu^s_U)(\Omega^\cdot(G))(U)\ .
\end{eqnarray*}

Since $\nu_U^s(\Omega^\cdot(G))$ is the de Rham complex of the manifold $U$ its higher cohomology sheaves vanish by the Poincar\'e Lemma. This implies that
$H^k(\Omega^\cdot(G))\cong 0$ for $k\ge 1$.

Furthermore, it is well-known that
$H^0(\nu_U(\Omega^\cdot(G)))\cong i^\sharp_{(U)}\R_{(U)}$.
It follows from the observation 
$$\nu_U^s\circ  i^\sharp_{\bG}\R_{\bG}\cong  i^\sharp_{(U)}\R_{(U)}$$ (proved by arguments similar as  above) that
$i^\sharp_{\bG}\R_{\bG}\cong H^0(\Omega^\cdot(G))$.
\hB

\subsubsection{}\label{klhu}

By Lemmas \ref{resd97} and \ref{red98} the complex $\Omega^\cdot(G)$ is a flabby resolution of $i^\sharp \R_\bG$. Therefore $$Ri(i^\sharp \R_\bG)\cong i(\Omega^\cdot(G))$$ in $D^+(\Pr_\Ab\bG)$.
By Lemma \ref{schluss} we have the isomorphism 
\begin{equation}\label{derham-model}(Rf_*\circ Ri)(i^\sharp\R_G)\cong C_A(\Omega^\cdot(G))
\end{equation}
 in $D^+(\Pr_{\Ab}\bX)$.
%Observe that $C_A(\Omega^\cdot(G))$ is actually a complex of sheaves.

\subsection{Calculation for $U(1)$-gerbes}

\subsubsection{}

In this subsection we specialize the situation of \ref{gensi} to the case where
$f:G\rightarrow X$ is a gerbe with band $U(1)$ according to the definition in \ref{gerbw} over a manifold $X$.
We thus can  assume that $A\to X$ is an atlas obtained from a covering of $X$ by open subsets such that the lift $s:A\rightarrow G$ is an atlas of $G$. The $U(1)$-central extension of groupoids in manifolds $(A\times_GA\Rightarrow A)\to (A\times_XA\Rightarrow A)$ (we forget the structure maps to $G$ for the moment) is the picture of a gerbe as presented in \cite{MR1876068}. In order to compare the sheaf theoretic construction of the cohomology of $G$  with the twisted de Rham complex we must choose some additional geometric structure on $G$, namely a connection in the sense of \cite{MR1876068}. The comparison map will depend on this choice.

\subsubsection{}

A connection on the  gerbe $f:G\rightarrow X$ consists of a pair $(\alpha,\beta)$,
where $\alpha\in \Omega^1(A\times_GA)$ is a connection one-form on  the $U(1)$-bundle  $A\times_GA\rightarrow A\times_XA$, and $\beta\in \Omega^2(A)$.
Observe that  $\Omega^2(A)$ and $\Omega^1(A\times_GA)$ are the first two spaces of the degree-two part of the graded commutative DG-algebra $C_A(\Omega^\cdot(G))(X)\cong \Omega^\cdot(G)(A^\cdot)$ discussed in \ref{alg} and \ref{localo}.  To be a connection the pair $(\alpha,\beta)$ is required to satisfy the following two conditions:
\begin{enumerate}
\item
$\delta \beta=d_{dR}\alpha$ (where $\delta$ is the \v{C}ech differential of the complex $\Omega^\cdot(A^\cdot)$, and $d_{dR}$ is the de Rham differential)
and \item $\delta\alpha=0$.
\end{enumerate} 
Note that $\delta d_{dR}\beta=0$ so that there is a unique $\lambda\in \Omega^3(X)$ which restricts to $d_{dR}\beta$. We have $d_{dR}\lambda=0$, and the class
$[\lambda]\in H^3(X;\R)$ represents the image under $H^3(X;\Z)\rightarrow H^3(X;\R)$ of the Dixmier-Douady class of the gerbe $G\rightarrow X$ (see \cite{MR1876068} for this fact and the existence of connections).

\subsubsection{} 

Let us choose a connection $(\alpha,\beta)$, and let $\lambda\in \Omega^3(X)$ be the associated closed three form. We consider $(\alpha,\beta)\in C_A(\Omega^\cdot(G))^2(X)$.

We consider the sheaf of complexes
$\Omega^\cdot[[z]]_\lambda$ on $\bX$ which associates to $(U\stackrel{i}{\rightarrow} X)\in \bX$ the complex
$$\Omega^\cdot(U)[[z]], \quad d_\lambda:=d_{dR}+ \lambda T\ ,$$
where $T:=\frac{d}{dz}$, $z$ has degree two, and $\lambda$ acts by right multiplication by $i^*\lambda$. In particular we have $d_\lambda z=i^*\lambda$. 
Note that $\Omega^\cdot[[z]]_\lambda$ is a sheaf of left $\Omega_X^\cdot$-$DG$-algebras.

\subsubsection{}\label{free}
% We have a global section $z_X\in \Omega^\cdot[[z]]_\lambda(X)$.
Observe that $z\in \Omega^\cdot[[z]]_\lambda(X)$ is central.
Let $L\in \Pr_{\Ab}\bX$ be a presheaf of graded unital central  $\Omega_X^\cdot$-algebras. A map of presheaves of graded unital central $\Omega_X^\cdot$-algebras $\phi:i\Omega^\cdot[[z]]_\lambda\rightarrow L$ determines a section $\phi(z)\in L(X)$. Vice versa, given
a section $ l\in L(X)$ of degree two, there is a unique
map of presheaves of graded unital central $\Omega_X^\cdot$-algebras $\phi:i\Omega^\cdot[[z]]_\lambda\rightarrow L$ such that $\phi(z)=l$. For $(U\stackrel{i}{\to} X)\in \bX$ the map $\phi_U:i\Omega^\cdot[[z]]_\lambda(U)\to L(U)$ is given by
$$\phi_U(\sum_{k\ge 0} \omega_k z^k):=\sum_{k\ge 0} \omega_k i^*(l)^k\ ,$$
where $i^*:L(X)\to L(U)$ is determined by the presheaf structure of $L$.

If $(L,d^L)$ is a presheaf of $DG$-algebras over  $\Omega_X^\cdot$,
then $\phi$ is a homomorphism of $DG$-algebras over  $\Omega_X^\cdot$ if and only if
$d^L l = \lambda$.

\subsubsection{}\label{prode}

\begin{prop}\label{er3}
We have an isomorphism
$$\Omega^\cdot[[z]]_\lambda\cong i^\sharp C_A(\Omega^\cdot(G))$$
in $D^+(\Sh_{\Ab}\bX)$.
\end{prop}
\proof
$C_A(\Omega^\cdot(G))$ is a presheaf of $DG$-algebras by the \ref{shdg}.
Given $(U\rightarrow X)\in \bX$ we have a natural projection $\pi:A_U^0\to U$
(see \ref{localo} for the notation). It induces a homomorphism of $DG$-algebras 
$\Omega_X^{\cdot}(U)\to (\ker(\delta):\Omega^\cdot(G)(A_U^0)\to \Omega^\cdot(G)(A_U^1))$ and therefore on 
$\Omega^\cdot(G)(A_U^\cdot)$ the structure of an $\Omega_X^\cdot(U)$-$DG$-module (see \ref{shdg}). In this way $C_A(\Omega^\cdot(G))$ becomes a sheaf of central $\Omega_X^\cdot$-$DG$-algebras.

By the discussion in \ref{free} we can
define a map of presheaves of central $\Omega_X^\cdot$-algebras
$$\tilde \phi:i\Omega^\cdot[[z]]_\lambda\rightarrow C_A(\Omega^\cdot(G))$$ such that
$\tilde \phi(z)=(\alpha,\beta)\in C_A(\Omega^\cdot(G))^2(X)$.
Because of $d \tilde \phi(z)=d(\alpha,\beta)=\lambda$, the map $\tilde \phi$ is a map
of presheaves of $DG$-algebras over $\Omega_X^\cdot$, hence in particular a map of presheaves of complexes.

We let
$$\phi:\Omega^\cdot[[z]]_\lambda\stackrel{\sim}{\to}i^\sharp \circ i \Omega^\cdot[[z]]_\lambda \stackrel{i^\sharp \tilde \phi}{\to } i^\sharp C_A(\Omega^\cdot(G))$$
be the induced map, where the first isomorphism exists since $\Omega^\cdot[[z]]_\lambda$ is a complex of sheaves.

\subsubsection{}
It remains to show that $\phi$ is a quasi-isomorphism of complexes of sheaves.
This can be shown locally. We can therefore assume that $X$ is contractible. We then  have a pull-back diagram
$$\xymatrix{G\ar[r]^q\ar[d]^f&[*/S^1]\ar[d]^g\\X\ar[r]^p&\mbox{*}}\ .$$
Since $p$ is smooth, so is $q$.
By Lemma \ref{derpulb} we have a canonical isomorphism
\begin{equation}\label{yxcv12}p^*\circ Rg_*\circ Ri\stackrel{\sim}{\to} Rf_*\circ q^*\circ Ri\ .\end{equation}

Applying (\ref{yxcv12}) to $i^\sharp \R_{\Site([*/S^1])}$ we obtain
\begin{equation}\label{waqw1}
p^*\circ Rg_*\circ  Ri(i^\sharp \R_{\Site([*/S^1])})\stackrel{\sim}{\rightarrow} Rf_*\circ q^* \circ Ri(i^\sharp\R_{\Site([* /S^1])})
\end{equation}
in $D^+(\Pr_\Ab\bX)$.
We now use (see \ref{klhu}) that   $$Ri(i^\sharp\R_{\Site([* /S^1])})\cong i(\Omega^\cdot{([*/S^1])})\ .$$
By the calculation of $q^*$ in Lemma \ref{poi451} and the definition of the de Rham complex we  have
$$q^*\circ i(\Omega^\cdot{([*/S^1])})\cong i(\Omega^\cdot(G))\ .$$
Therefore in $D^+(\Pr_\Ab\bG)$ we have
$$q^*\circ Ri(i^\sharp\R_{\Site([* /S^1])})\cong  i(\Omega^\cdot(G)) \stackrel{\ref{klhu}}{\cong} Ri(i^\sharp\R_\bG)\ .$$
It follows by \ref{klhu} that
\begin{equation}
\label{waqw2}Rf_*\circ q^* \circ Ri(i^\sharp\R_{\Site([* /S^1])})\cong Rf_*\circ Ri(i^\sharp\R_\bG)\cong C_A(\Omega^\cdot(G))\ .
\end{equation}

\subsubsection{}

We now must calculate the cohomology of the gerbe $[*/S^1]$ with real coefficients.
\begin{lem}\label{bas13}
We have an isomorphism 
$$i^\sharp\circ Rg_*\circ Ri(i^\sharp\R_{\Site([*/S^1])})\cong i^\sharp (\R[[z]]_{\Site(*)})\ ,$$ where $z$ has degree two.
\end{lem}
\proof
We choose the atlas $A:=*\to [*/S^1]$ and use the isomorphism
$$Rg_*\circ Ri(i^\sharp\R_{\Site([*/S^1])})\cong C_A(\Omega^\cdot([*/S^1]))\in D^+({\Pr}_\Ab\Site(*))\ .$$ Note that $\Site(*)$ is the category of smooth manifolds.
Let $U$ be a smooth manifold. We have  $$A_U\cong U\to U\times [*/S^1]\cong [*/S^1]_U $$and 
$$C_A(\Omega^\cdot([*/S^1]))(U)\cong \Omega^\cdot(A_U^\cdot)\ ,$$
where $$A_U^p=\underbrace{A_U\times_{[*/S^1]_U }\dots\times _{[*/S^1]_U}A_U}_{p+1 \:factors}\cong U\times (\underbrace{*\times_{[*/S^1] }\dots\times _{[*/S^1]}*}_{p+1 \:factors}) \cong U\times (S^1)^{p}\ .$$
The simplicial manifold
$$A_U^\cdot\cong U\times (S^1)^\cdot$$ is the simplicial model of the space
$U\times BS^1$, where $BS^1$ is the classifying space of the group $S^1$.
We can use the simplicial de Rham complex in order to calculate its cohomology.
Note that $H^*(BS^1,\R)\cong \R[[z]]$ with $z$ in degree two.
Let us fix a form $\zeta\in (\Omega^\cdot ((S^1)^\cdot)_{tot}^2$ which represents the generator $z$. Then we define a map
$$\mu_U: \Omega^\cdot(U)[[z]]\to \Omega^\cdot(U\times (S^1)^\cdot)$$
by $$\mu(\omega z^k):=\omega\wedge \zeta^k\ .$$
This map induces a quasi-isomorphism of complexes of abelian groups.
The family of maps $\mu_U$ for varying $U$ defines a quasi-isomorphism of complexes of presheaves $\mu:i\Omega(*)[[z]]\to C_A(\Omega^\cdot([*/S^1]))$.
It induces the quasi-isomorphism of complexes of sheaves
$$\Omega(*)[[z]]\cong i^\sharp\circ i \Omega(*)[[z]]\stackrel{i^\sharp \mu}{\cong} i^\sharp C_A(\Omega^\cdot([*/S^1]))\ .$$
Finally observe that the canonical map
$$i^\sharp \R[[z]]_{\Site(*)}\to\Omega(*)[[z]]$$
is a quasi-isomorphism by Lemma \ref{red98} 
\hB

\subsubsection{}

 It follows from Lemma \ref{bas13} by applying $p^*\circ Ri$ that
$$p^*\circ Ri\circ i^\sharp \circ Rg_* \circ Ri(i^\sharp\R_{\Site([*/S^1])})\cong  p^*\circ Ri\circ  i^\sharp\R[[z]]_{\Site(*)}\ .$$  We now use the second assertion of Lemma \ref{somecomu} in order to commute $Ri\circ i^\sharp$ with $p^*$. We get
$$Ri\circ i^\sharp \circ p^*\circ  Rg_* \circ Ri(i^\sharp\R_{\Site([*/S^1])})\cong   Ri\circ  i^\sharp\circ p^* \R[[z]]_{\Site(*)}\ .$$
We now apply $i^\sharp$ and use that $i^\sharp \circ Ri\cong \id$ in order to drop the functor $Ri$ and get the quasi-isomorphism
$$i^\sharp \circ p^*\circ  Rg_* \circ Ri(i^\sharp\R_{\Site([*/S^1])})\cong    i^\sharp\circ p^* \R[[z]]_{\Site(*)}\ .$$

By the explicit description of $p^*$ given in the proof of Lemma \ref{pi} we see that
$$ p^*(\R[[z]]_{\Site(*)})\cong  \R[[z]]_\bX\ .$$
We thus have a quasi-isomorphism
\begin{equation}\label{waqw3}
i^\sharp \circ p^*\circ Rg_* \circ Ri(i^\sharp\R_{\Site([*/S^1])})\cong    i^\sharp(\R[[z]]_\bX)\ .
\end{equation}

Combining the isomorphisms (\ref{waqw1}), (\ref{waqw2}) and  (\ref{waqw3}) we obtain a quasi-isomorphism
$$i^\sharp C_A(\Omega^\cdot(G))\cong i^\sharp(\R[[z]]_\bX)\ .$$
In particular we see that $z$ generates the cohomology.

\subsubsection{}

Since $X$ is contractible we find 
$\gamma\in \Omega^2(X)$ such that $d_{dR}\gamma=\lambda$.
We define a map of complexes of sheaves
$$\psi:i^\sharp\R[[z]]_\bX\rightarrow \Omega^\cdot[[z]]_0\stackrel{e^{-\gamma T}}{\rightarrow}  \Omega^\cdot[[z]]_\lambda\ .$$
The first map is given by the inclusion $i^\sharp\R_\bX\rightarrow \Omega_X^\cdot$ 
and is a quasi-isomorphism. The second map is an isomorphism of sheaves of complexes. Therefore $\psi$ is a quasi-isomorphism. Note that $\psi$ is multiplicative and $\psi(z)=z-\gamma$.
We further define $\kappa:i^\sharp\R[[z]]_\bX\rightarrow i^\sharp C_A(\Omega^\cdot(G))$ such that
$\kappa(z)=(\alpha,\beta-\gamma)=\phi(z-\gamma)$.
Then we have a commutative diagram
$$\xymatrix{ i^\sharp\R[[z]]_\bX\ar[rr]^\kappa\ar[dr]^\psi&&i^\sharp C_A(\Omega^\cdot(G))\\&\Omega^\cdot[[z]]_\lambda\ar[ur]^\phi}\ .$$ If we show that $\kappa$ is a quasi-isomorphism, then since
$\psi$ is a quasi-isomorphism, $\phi$ must be a quasi-isomorphism, too.
It suffices to see that $\kappa(z):=(\alpha,\beta-\gamma)$ represents a non-trivial cohomology class. Assume that it is a boundary locally on $(U\rightarrow X)\in \bX$. Then there exists
$x\in \Omega^0(G)(A_U^1)$ and $y\in \Omega^1(G)(A_U^0)$ such that
$\delta x=0$, $d_{dR}x+\delta y=\alpha_U$ and $d_{dR} y=(\beta-\gamma)_U$
(the subscript indicates that the forms are pulled back to $A_U^*$).
By exactness of the $\delta$-complex we can in fact assume that $x=0$.
But then the equation $\delta y=\alpha_U$ is impossible since
$\delta y$ vanishes on vertical vectors on the bundles $A_U^1\rightarrow A_U^0$ given by the source and range projections while $\alpha$ as a connection form is non-trivial on those vectors.\hB

\subsubsection{}\label{fini}

We now finish the proof of Theorem \ref{main}.
We combine Proposition \ref{er3} with \ref{klhu} and the fact that $f_*$ preserves sheaves (Lemma \ref{shpre}) in order to get 
$$ i^\sharp\circ Rf_*\circ Ri(i^\sharp \R_\bG)\cong i^\sharp C_A(\Omega^\cdot(G))\cong \Omega^\cdot[[z]]_\lambda\ .$$

%\section{The Chern characer}

%\subsection{The lifting Gerbe}

%\subsubsection{}

%In the present Subsection we describe a canonical construction which associates to a $PU$-principal bundle $P\rightarrow X$ over a space $X$ a topological gerbe $G(P)\rightarrow X$ with band $U(1)$. It will be called the lifting gerbe.

%\subsubsection{}

%We work in the category $\Top$ and will consider stacks in topological spaces.
%We choose some separable infinite dimensional complex Hilbert space $H$ and let $U$ be the topological group of unitary automorphisms of $H$ with the topology induced by the operator norm. The group $U(1)$ acts on $H$ by complex multiplication. It maps isomorphically to the center of $U$. We define $PU:=U/U(1)$, the projective unitary group of $H$. It is a topological group, and we can consider $PU$-principal bundles $P\rightarrow X$.

%\subsubsection{}

%Let $P\to X$ be a $PU$-principal bundle. Observe that $U$ acts on $P$ via the projection $U\to PU$.
%We define a presheaf $G:=G(P)$ of groupoids on $\Top$ as follows.
%For a space $T\in \Top$ the objects of $G(T)$ are the diagrams
%$$\xymatrix{Q\ar[r]\ar[d]&P\ar[d]\\
%T\ar[r]&X}\ ,$$ where $Q\to T$ is a $U$-principal bundl and
%$Q\to P$ is $U$-equivariant.

%A morphism between two such objects
%$$\xymatrix{Q\ar[r]\ar[d]&P\ar[d]\\
%T\ar[r]&X}\ ,\xymatrix{Q^\prime \ar[r]\ar[d]&P\ar[d]\\
%T\ar[r]&X}\ ,$$
% is an isomorphism of $U$-principal bundles
%$Q\to Q^\prime$ over $T$ which is compatible with the maps to $P$.

\hB


\begin{thebibliography}{10}

\bibitem{math.KT/0510674}
Michael Atiyah and Graeme Segal.
\newblock {Twisted K-theory and cohomology}.

\bibitem{MR2172633}
Michael Atiyah and Graeme Segal.
\newblock Twisted {$K$}-theory.
\newblock {\em Ukr. Mat. Visn.}, 1(3):287--330, 2004.

\bibitem{MR2172499}
K.~Behrend.
\newblock Cohomology of stacks.
\newblock In {\em Intersection theory and moduli}, ICTP Lect. Notes, XIX, pages
  249--294 (electronic). Abdus Salam Int. Cent. Theoret. Phys., Trieste, 2004.

\bibitem{math.DG/0605694}
Kai Behrend and Ping Xu.
\newblock {Differentiable Stacks and Gerbes}.

\bibitem{MR2183389}
Kai~A. Behrend.
\newblock On the de {R}ham cohomology of differential and algebraic stacks.
\newblock {\em Adv. Math.}, 198(2):583--622, 2005.

\bibitem{MR1911247}
Peter Bouwknegt, Alan~L. Carey, Varghese Mathai, Michael~K. Murray, and Danny
  Stevenson.
\newblock Twisted {$K$}-theory and {$K$}-theory of bundle gerbes.
\newblock {\em Comm. Math. Phys.}, 228(1):17--45, 2002.

\bibitem{MR2080959}
Peter Bouwknegt, Jarah Evslin, and Varghese Mathai.
\newblock {$T$}-duality: topology change from {$H$}-flux.
\newblock {\em Comm. Math. Phys.}, 249(2):383--415, 2004.

\bibitem{MR1197353}
Jean-Luc Brylinski.
\newblock {\em Loop spaces, characteristic classes and geometric quantization},
  volume 107 of {\em Progress in Mathematics}.
\newblock Birkh\"auser Boston Inc., Boston, MA, 1993.

% \bibitem{bssf}
% U.~Bunke, Th. Schick, and M.~Spitzweck.
% \newblock {Foundations of sheaf theory on topological stacks}.
% \newblock In preparation.

\bibitem{bss1}
U.~Bunke, Th. Schick, and M.~Spitzweck.
\newblock {Inertia and delocalized twisted cohomology}.
\newblock In preparation.

\bibitem{bssm}
U.~Bunke, Th. Schick, and M.~Spitzweck.
\newblock {$T$-duality and periodic twisted cohomology}.
\newblock In preparation.

\bibitem{math.AT/0206257}
Daniel~S. Freed, Michael~J. Hopkins, and Constantin Teleman.
\newblock {Twisted equivariant K-theory with complex coefficients}.

\bibitem{heinloth}
Jochen Heinloth.
\newblock Survey on topological and smooth stacks.
\newblock In {\em Mathematisches Institut G{\"o}ttingen, WS04-05 (Y. Tschinkel,
  ed.)}, pages 1--31. 2005.

\bibitem{MR1876068}
Nigel Hitchin.
\newblock Lectures on special {L}agrangian submanifolds.
\newblock In {\em Winter School on Mirror Symmetry, Vector Bundles and
  Lagrangian Submanifolds (Cambridge, MA, 1999)}, volume~23 of {\em AMS/IP
  Stud. Adv. Math.}, pages 151--182. Amer. Math. Soc., Providence, RI, 2001.

\bibitem{MR2122155}
Michael Joachim.
\newblock Higher coherences for equivariant {$K$}-theory.
\newblock In {\em Structured ring spectra}, volume 315 of {\em London Math.
  Soc. Lecture Note Ser.}, pages 87--114. Cambridge Univ. Press, Cambridge,
  2004.

\bibitem{MR1299726}
Masaki Kashiwara and Pierre Schapira.
\newblock {\em Sheaves on manifolds}, volume 292 of {\em Grundlehren der
  Mathematischen Wissenschaften [Fundamental Principles of Mathematical
  Sciences]}.
\newblock Springer-Verlag, Berlin, 1994.
\newblock With a chapter in French by Christian Houzel, Corrected reprint of
  the 1990 original.

\bibitem{MR1771927}
G{\'e}rard Laumon and Laurent Moret-Bailly.
\newblock {\em Champs alg\'ebriques}, volume~39 of {\em Ergebnisse der
  Mathematik und ihrer Grenzgebiete. 3. Folge. A Series of Modern Surveys in
  Mathematics}.
\newblock Springer-Verlag, Berlin, 2000.

\bibitem{math.KT/0404329}
Varghese Mathai and Danny Stevenson.
\newblock {On a generalized Connes-Hochschild-Kostant-Rosenberg theorem}.
\newblock Advances in Mathematics, vol. 200 no. 2 (2006) 1-33.

\bibitem{MR1977885}
Varghese Mathai and Danny Stevenson.
\newblock Chern character in twisted {$K$}-theory: equivariant and holomorphic
  cases.
\newblock {\em Comm. Math. Phys.}, 236(1):161--186, 2003.

\bibitem{math.AT/0411656}
J.~P. May and J.~Sigurdsson.
\newblock {Parametrized homotopy theory}.

\bibitem{MR0494077}
J.~Peter May.
\newblock {\em {$E\sb{\infty }$} ring spaces and {$E\sb{\infty }$} ring
  spectra}.
\newblock Springer-Verlag, Berlin, 1977.
\newblock With contributions by Frank Quinn, Nigel Ray, and J\o rgen Tornehave,
  Lecture Notes in Mathematics, Vol. 577.

\bibitem{math.DG/0306176}
David Metzler.
\newblock {Topological and Smooth Stacks}.

\bibitem{MR1405064}
M.~K. Murray.
\newblock Bundle gerbes.
\newblock {\em J. London Math. Soc. (2)}, 54(2):403--416, 1996.

\bibitem{MR1794295}
Michael~K. Murray and Daniel Stevenson.
\newblock Bundle gerbes: stable isomorphism and local theory.
\newblock {\em J. London Math. Soc. (2)}, 62(3):925--937, 2000.

\bibitem{math.AG/0503247}
Behrang Noohi.
\newblock {Foundations of Topological Stacks I}.

\bibitem{olsson}
Martin Olsson.
\newblock Sheaves on artin stacks.

\bibitem{MR1401424}
Dorette~A. Pronk.
\newblock Etendues and stacks as bicategories of fractions.
\newblock {\em Compositio Math.}, 102(3):243--303, 1996.

\bibitem{MR1317816}
G{\"u}nter Tamme.
\newblock {\em Introduction to \'etale cohomology}.
\newblock Universitext. Springer-Verlag, Berlin, 1994.
\newblock Translated from the German by Manfred Kolster.

\bibitem{math.KT/0306138}
Jean-Louis Tu, Ping Xu, and Camille Laurent-Gengoux.
\newblock Twisted k-theory of differentiable stacks.

\end{thebibliography}
\end{document}